\documentclass{gtart_h}

\def\ifplaintex{\expandafter\ifx\csname documentclass\endcsname\relax}

\def\gtp{{\mathsurround=0pt\it $\cal G\mskip-2mu$eometry \&\ 
$\cal T\!\!$opology $\cal P\!$ublications}}  

\def\recd{{\small Received:\qua\receiveddate\ifx\reviseddate\relax
\else\qquad Revised:\qua\reviseddate\fi\par}} 

\def\lognumber#1{\def\thelognumber{#1}}
\def\volumenumber#1{\def\thevolumenumber{#1}}
\def\volumeyear#1{\def\thevolumeyear{#1}}
\def\papernumber#1{\def\thepapernumber{#1}}
\def\pagenumbers#1#2{\def\startpage{#1}\def\finishpage{#2}}
\def\published#1{\def\publishdate{#1}}

\def\received#1{\def\receiveddate{#1}}
\def\revised#1{\def\reviseddate{#1}}
\def\accepted#1{\def\accepteddate{#1}}
\def\asciititle#1{\def\theasciititle{#1}}

\long\def\asciiabstract#1{\long\def\theasciiabstract{#1}}
\def\asciikeywords#1{\def\theasciikeywords{#1}}

\let\\\par\let\thelognumber\relax\let\thevolumenumber\relax
\let\thepapernumber\relax\let\thevolumeyear\relax\let\startpage\relax
\let\finishpage\relax\let\publishdate\relax\let\receiveddate\relax
\let\reviseddate\relax\let\accepteddate\relax\let\theasciititle\relax
\let\theasciiauthors\relax
\let\theasciiabstract\relax\let\theasciikeywords\relax

\let\theasciiemail\relax

\ifplaintex
\font\logobig=cmssbx10 scaled 3836
\font\logomed=cmssbx10 scaled 2557
\else
\font\logobig=cmssbx10 scaled 4200
\font\logomed=cmssbx10 scaled 2800
\fi

\long\def\makeagttitle{   
\count0=\startpage
\agt\hfill      
\hbox to 45truept{\vbox to 0pt{\vglue -13truept{\logomed A\kern -.37em{\logobig 
T}\kern -.38em G}\vss}\hss}
\break
{\small Volume \thevolumenumber\ (\thevolumeyear)
\startpage--\finishpage\nl
Published: \publishdate}

\vglue .25truein

{\parskip=0pt\leftskip 0pt plus
1fil\def\\{\par\smallskip}{\Large\bf\thetitle}\par\medskip} \vglue
0.05truein

{\parskip=0pt\leftskip 0pt plus 1fil\def\\{\par}{\sc\theauthors}
\par\medskip}%
 
\vglue 0.03truein 

{\small\leftskip 25truept\rightskip 25truept{\bf Abstract}\stdspace\theabstract

{\bf AMS Classification}\stdspace\theprimaryclass
\ifx\thesecondaryclass\relax\else; \thesecondaryclass\fi\par
{\bf Keywords}\stdspace \thekeywords\par}\vglue 7truept

}   

\ifplaintex
\font\phead=cmsl9 scaled 950
\font\pnum=cmbx10 scaled 913
\font\pfoot=cmsl9 scaled 950
\headline{\vbox to 0pt{\vskip -4.5mm\line{\small\phead\ifnum
\count0=\startpage ISSN 1472-2739 (on-line) 1472-2747 (printed)
\hfill {\pnum\folio}\else\ifodd\count0\def\\{ }%
\ifx\theshorttitle\relax\thetitle\else\theshorttitle\fi\hfill{\pnum\folio}
\else\def\\{ and }{\pnum\folio}\hfill\ifx\theshortauthors\relax\theauthors
\else\theshortauthors\fi\fi\fi}\vss}}
\footline{\vbox to 0pt{\vglue 0mm\line{\small\pfoot\ifnum\count0=\startpage
\copyright\ \gtp\hfill\else
\agt, Volume \thevolumenumber\ (\thevolumeyear)\hfill\fi}\vss}}
\else
\font=cmsl9 scaled 1050
\font\lnum=cmbx10 
\font=cmsl9 scaled 1050
\makeatletter
\def\@oddhead{{\small\lhead\ifnum\count0=\startpage ISSN 1472-2739 
(on-line) 1472-2747 (printed)\hfill {\lnum\number\count0}\else\ifodd\count0
\def\\{ }\ifx\theshorttitle\relax \thetitle \else\theshorttitle\fi\hfill
{\lnum\number\count0}\else\def\\{ and }{\lnum\number\count0}
\hfill\ifx\theshortauthors\relax 
\theauthors\else\theshortauthors\fi\fi\fi}}\def\@evenhead{\@oddhead}
\def\@oddfoot{\small\lfoot\ifnum\count0=\startpage\copyright\ \gtp\hfill\else
\agt, Volume \thevolumenumber\ (\thevolumeyear)\hfill\fi}
\def\@evenfoot{\@oddfoot}
\makeatother
\fi
\let\maketitlepage\makeagttitle
\let\makeshorttitle\maketitlepage
\let\maketitle\maketitlepage

\newwrite\gtoutfile
\long\gdef\makeheadfile{  
{\def\\{, }\def\s{ }
\immediate\openout\gtoutfile head.xxx
\immediate\write\gtoutfile{Proxy-for: \ifx\theasciiauthors\relax
\theauthors\else\theasciiauthors\fi\s<\ifx\theasciiemail\relax\theemail\else\theasciiemail\fi>}
\immediate\write\gtoutfile{\noexpand\\}
\immediate\write\gtoutfile{Authors: \ifx\theasciiauthors\relax
\theauthors\else\theasciiauthors\fi}
{\def\\{ }\immediate\write\gtoutfile{Title: \ifx\theasciititle\relax
\thetitle\else\theasciititle\fi}}
\immediate\write\gtoutfile{Subj-class: GT or SG, GR etc}
\immediate\write\gtoutfile{MSC-class: \theprimaryclass\ifx\thesecondaryclass\relax\else, \thesecondaryclass\fi}
\immediate\write\gtoutfile{Journal-ref: Algebr. Geom. Topol. \thevolumenumber\s
(\thevolumeyear) \startpage-\finishpage}
\immediate\write\gtoutfile{Comments: Published by Algebraic and
Geometric Topology at}
\immediate\write\gtoutfile{\s\s\s  http://www.maths.warwick.ac.uk/agt/AGTVol\thevolumenumber/agt-\thevolumenumber-\thepapernumber.abs.html}
\immediate\write\gtoutfile{\noexpand\\}
\immediate\write\gtoutfile{}
\ifx\theasciiabstract\relax
\immediate\write\gtoutfile{\theabstract}\else
\immediate\write\gtoutfile{\theasciiabstract}\fi
\immediate\write\gtoutfile{}
\immediate\write\gtoutfile{\noexpand\\}
\immediate\write\gtoutfile{}
\immediate\closeout\gtoutfile}}  

\def\maketitlepage{\makeagttitle\makeheadfile}
\let\makeshorttitle\maketitlepage
\let\maketitle\maketitlepage

\lognumber{23}
\volumenumber{4}
\volumeyear{2004}
\papernumber{23}
\published{8 July 2004}
\pagenumbers{473}{520}
\received{18 December 2003}
\revised{26 May 2004}
\accepted{6 July 2004}

\usepackage{amssymb}
\usepackage{amsmath}
\usepackage{graphicx}

\def\bu{_\bullet}

\def\xrarrow{\xrightarrow} 

\def\noteq{\neq}

\def\grad{\nabla}
\def\<{\left<}
\def\>{\right>}

\def\nequiv{\equiv\!\!\!\!\!/\ }

\DeclareMathOperator{\codim}{codim}

\DeclareMathOperator{\sgn}{sgn} \DeclareMathOperator{\val}{val}
\DeclareMathOperator{\Homeo}{Homeo}

\newcommand{\field}[1]{\mathbb{#1}}
\newcommand{\ZZ}{\ensuremath{{\field{Z}}}}
\newcommand{\CC}{\ensuremath{{\field{C}}}}
\newcommand{\RR}{\ensuremath{{\field{R}}}}
\newcommand{\QQ}{\ensuremath{{\field{Q}}}}
\newcommand{\NN}{\ensuremath{{\field{N}}}}

\def\l{\ell}
\def\ll{\lambda}
\def\LL{\Lambda}
\newcommand{\A}{\ensuremath{{\mathcal{A}}}}

\newcommand{\C}{\ensuremath{{\mathcal{C}}}}

\newcommand{\F}{\ensuremath{{\mathcal{F}}}}
\newcommand{\G}{\ensuremath{{\mathcal{G}}}}

\newcommand{\nL}{\ensuremath{{\mathcal{L}}}}

\newcommand{\N}{\ensuremath{{\mathcal{N}}}}
\newcommand{\nP}{\ensuremath{{\mathcal{P}}}}

\newcommand{\W}{\ensuremath{{\mathcal{W}}}}

\newcommand{\Z}{\ensuremath{{\mathcal{Z}}}}

\def\g{\gamma}
\def\d{\partial}

\def\e{\epsilon}
\def\f{\phi}
\def\k{\kappa}

\def\s{\sigma}
\def\Sig{\Sigma}
\def\t{\tau}

\def\w{\omega}

\def\FG{{\F\G}}
\def\Fat{{\F at}}
\def\halfbig{A_{2k+3}^{1/2}}
\def\halfsmall{A_{2k+2}^{1/2}}
\def\wt{\widetilde}

\def\st{\,|\,}

\def\Gam{\Gamma}

\def\lewq{\leq}

\newtheorem{thm}{Theorem}[section]
\newtheorem{lem}[thm]{Lemma}
\newtheorem{cor}[thm]{Corollary}
\newtheorem{prop}[thm]{Proposition}
\theoremstyle{definition}
\newtheorem{defn}[thm]{Definition}
\newtheorem{rem}[thm]{Remark}

\begin{document}
\title{Combinatorial Miller--Morita--Mumford classes\\and Witten cycles}
\asciititle{Combinatorial Miller-Morita-Mumford classes and Witten cycles}
\authors{Kiyoshi Igusa}                
\address{Department of Mathematics, Brandeis University\\Waltham, MA
02454-9110, USA}
\email{igusa@brandeis.edu}

\begin{abstract}  
We obtain a combinatorial formula for the Miller--Morita--Mum\-ford
classes for the mapping class group of punctured surfaces and
prove Witten's conjecture that they are proportional to the dual
to the Witten cycles. The proportionality constant is shown to be
exactly as conjectured by Arbarello and Cornalba
\cite{[Arbarello-Cornalba:96]}. We also verify their conjectured
formula for the leading coefficient of the polynomial expressing
the Kontsevich cycles in terms of the Miller--Morita--Mumford
classes.
\end{abstract}
\asciiabstract{%
We obtain a combinatorial formula for the Miller-Morita-Mumford
classes for the mapping class group of punctured surfaces and prove
Witten's conjecture that they are proportional to the dual to the
Witten cycles. The proportionality constant is shown to be exactly as
conjectured by Arbarello and Cornalba [J. Alg. Geom. 5 (1996)
705-749]. We also verify their conjectured formula for the leading
coefficient of the polynomial expressing the Kontsevich cycles in
terms of the Miller-Morita-Mumford classes.}



\asciikeywords{Mapping class group, fat graphs, ribbon graphs,
Miller-Morita-Mumford classes, tautological classes, Witten
conjecture, Stasheff associahedra}                    

\maketitle 

%
%
\section*{Introduction}
\addcontentsline{toc}{section}{Introduction}

The Miller--Morita--Mumford classes were defined by David Mumford
\cite{[Mumford83:MMM_class]} as even dimensional cohomology
classes on the Deligne--Mumford compactification of the moduli
space of Riemann surfaces of a fixed genus~$g$. These are also
referred to as \emph{tautological classes}.

Around the same time Shigeyuki Morita \cite{[Morita84]} defined
\emph{characteristic classes} for oriented surface bundles. These
are topologically defined integer cohomology classes for the
\emph{mapping class group} $M_g$, ie, the group of isotopy
classes of orientation preserving diffeomorphisms of a fixed
Riemann surface $\Sig_g$ of genus~$g$:
\[
    \k_k\in H^{2k}(M_g;\ZZ)
\]
The mapping class group has rational cohomology isomorphic to that
of the (uncompactified) moduli space. Ed Miller
\cite{[Miller86:MMM]} showed that Mumford's tautological classes
correspond under this isomorphism to these topologically defined
classes. Using work of John Harer \cite{[Harer85:Stability]}
Miller showed that these cohomology classes $\k_k$ are
algebraically independent in the stable range. This result was
also obtained independently by Morita \cite{[Morita87]}.

In this paper we take the topological viewpoint. We note that the
topological and algebraic geometric definitions of the
Miller--Morita--Mumford classes agree up to a sign of $(-1)^{k+1}$.
(See \cite{[MoritaSurvey]}.)

The mapping class group $M_g^s$ of genus $g$ surfaces with
$s\geq1$ boundary components (which we are allowed to rotate and
permute) is classified by a space of \emph{fat graphs} (also
called \emph{ribbon graphs}). These are defined to be finite
graphs where all vertices have valence 3 or more together with a
cyclic ordering of the half edges incident to each vertex. E.
Witten conjectured that the Miller--Morita--Mumford classes were
dual to certain $2k$--cycles in the space of fat graphs. We call
these the \emph{Witten cycles} and denote them by $W_{k}$. (See \ref{subsec:Witten cycle}.) M.
Kontsevich in \cite{[Kontsevich:Airy]} constructed other cycles in
the space of fat graphs and conjectured that they could be
expressed in terms of the Miller--Morita--Mumford classes.

Robert Penner \cite{[Penner91:PoincareDual]} verified Witten's
conjecture when $k=1$. However his calculation was off by a factor
of $2$. The correct statement for $k=1$ was given by E. Arbarello
and M. Cornalba \cite{[Arbarello-Cornalba:96]}:
\[
    \wt{\k}_1=\frac1{12}[W_1]^\ast
\]
where $\wt{\k}_1\in H^2(M_g^s;\ZZ)$ is the
\emph{adjusted Miller--Morita--Mumford class} defined in \ref{subsec:Morse}.

In this paper we prove the Witten conjecture for all $k\geq0$:

\begin{thm}\label{main theorem} The adjusted Miller--Morita--Mumford class $\wt{\k}_k$ is
related to the duals $[W_{k}]^\ast$ of the Witten cycles as
elements of $H^{2k}(M_g^s;\QQ)$ for all $k\geq0$ by
\[
    \wt{\k}_k=(-1)^{k+1}\frac{(k+1)!}{(2k+2)!}[W_{k}]^\ast
\]
as conjectured in \cite{[Arbarello-Cornalba:96]}.
\end{thm}

To prove this we construct an elementary combinatorial cocycle
representing the class $\wt{\k}_k$ (Theorem \ref{FG:thm:bk is proportional to Wk} ) and evaluate it on the cocycle $[W_{k}]$ (Theorem \ref{Calc:thm:proportionality constant between MMM and b}).

As an easy consequence of the above theorem (combining Theorem \ref{Calc:thm:proportionality constant between MMM and b} with Corollary \ref{FG:cor:Kontsevich cycles in terms of bks}) we obtain the
following.

\begin{cor}\label{intro:cor:Kontsevich cycles in terms of
MMM} The Kontsevich cycles $W_{k_1^{n_1}\cdots k_r^{n_r}}$
\cite{[Kontsevich:Airy]} are dual to polynomials in the adjusted
Miller--Morita--Mumford classes with leading terms as conjectured in
\cite{[Arbarello-Cornalba:96]}:
\[
    \left[
W_{k_1^{n_1}\cdots k_r^{n_r}}
    \right]^\ast
    =
    \prod_{i=1}^r\frac1{
n_i!}\left( 2\frac{(2k_i+1)!}{(-1)^{k_i+1}k_i!}
    \right)^{n_i}
    (\wt{\k}_{k_1})^{n_1}\cdots(\wt{\k}_{k_r})^{n_r}
    + \text{lower terms}.
\]
\end{cor}

Outline of the paper: In the first two sections we construct a $2k$ cocycle $c_\Z^k$ on the category of cyclically ordered set $\Z$ and show that it represents the $k^{\rm th}$ power of the Euler class on $|\Z|\simeq\CC P^\infty$. In Section 3 we use this to define a $2k$ cocycle $c_\Fat^k$ on the category of fat graphs $\Fat$ and show that it is proportional to the dual $[W_k]^\ast$ of the Witten cycle $[W_k]$. The  proportionality constant is computed at the end of this section using Stasheff associahedra. Section 4 uses Morse theory on surfaces associated to fat graphs to show that $[W_k]^\ast$ is proportional to the adjusted Miller--Morita--Mumford class $\wt{\k}_k$ on $|\Fat|\simeq \coprod BM_g^s$. The final section computes the proportionality constant between $\wt{\k}_k$ and $c_\Fat^k$ giving the main theorems as stated above.

{I would like to thank Robert Penner for explaining the
Witten conjecture to me in detail many years ago. More recently I
owe thanks to Jack Morava and Dieter Kotschick whose questions and
comments lead to this current successful attempt at this
conjecture. I am very grateful to Daniel Ruberman and Harry
Tamvakis for their help in developing the original ideas for the
combinatorial Miller--Morita--Mumford classes. Finally, I would like
to thank Nariya Kawazumi, Shigeyuki Morita and Kenji Fukaya and
everyone in Tokyo and Kyoto who helped and encouraged me during
the time that I was improving the results of this paper.

The main theorem of this paper has corollaries which are explained
in two other papers, the third jointly with Michael Kleber. These
subsequent papers (\cite{[I:GraphCoh]},
\cite{[IKleber:IncreasingTrees]}) also give formulas for the
``lower terms'' in the above expression. At the same time a paper
by Mondello \cite{[Mondello]} has appeared proving the same thing.
}
{This paper was written under NSF Grant DMS-0204386 and revised under DMS-0309480.}
%
%

\section{Cyclically ordered sets}

In this section we consider the category $\Z$ of cyclically
ordered sets and cyclically ordered monomorphisms. This category
is one of a number of well-known models for $\CC P^\infty$ and
therefore its integer cohomology is a polynomial algebra in its
Euler class
\[
    e_\Z\in H^2(|\Z|;\ZZ).
\]
We give an explicit rational cocycle $c_\Z^k$ for the $k^{\rm th}$ power
$e_\Z^k$ of this Euler class. The crucial point is that, in order
for $c_\Z^k$ to be non-zero on a $2k$--simplex
\[
    C_0\rightarrowtail C_1\rightarrowtail\cdots\rightarrowtail
    C_{2k}\ ,
\]
the cyclic sets $C_i$ must strictly increase in size:
\[
    |C_0|<|C_1|<\cdots<|C_{2k}|.
\]

\subsection{The category $\Z$ of cyclically ordered sets}

By a \emph{cyclically ordered set} we mean a finite non-empty set
$C$ with, say, $n$ elements, together with a cyclic permutation
$\s$ of $C$ of order $n$. Thus $C$ has $(n-1)!$ cyclic orderings.
We also use square brackets to denote cyclically ordered sets:
\[
    (C,\s)=[x,\s(x),\s^2(x),\cdots,\s^{n-1}(x)].
\]
To avoid set theoretic problems we assume that $C$ is a subset of
some fixed infinite set.

To each cyclically ordered set $(C,\s)$ we can associate an
oriented graph $S^1(C,\s)$ with one vertex for each element of $C$
and one directed edge $x\to y$ if $y=\s(x)$. Then $S^1(C,\s)$ is
homeomorphic to a circle.

Any monomorphism of cyclically ordered sets $f\co(C,\s)\to (D,\t)$
has a \emph{degree} given by
\[
    \deg(f)=\frac1{|D|}\sum_{x\in C}k(x)
\]
where $k=k(x)$ is the smallest positive integer so that
$f(\s(x))=\t^k(f(x))$. The degree of $f$ is also the degree of the
induced mapping $f_\ast\co S^1(C,\s)\to S^1(D,\t)$. However, note
that this is not a functor ($(fg)_\ast\noteq f_\ast g_\ast$ and
$\deg(fg)\noteq \deg(f)\deg(g)$ in general). For example, the
degree of any monomorphism $f$ is equal to $1$ if the domain has
$\leq 2$ elements.

Let $\Z$ be the category of cyclically ordered sets $(S,\s)$ with
morphisms $f\co(S,\s)$ $\to (T,\t)$ defined to be set monomorphisms
$f\co S\to T$ of degree $1$.

\subsection{Linearly ordered sets}

Let $\nL$ denote the category of finite, non-empty, linearly ordered
sets and order preserving monomorphisms. Then we have a functor
$J\co\nL\to\Z$ sending $C=(x_1<\cdots<x_n)$ to $(C,\s)$ where the
cyclic ordering $\s$ is given by $\s(x_n)=x_1$ and
$\s(x_i)=x_{i+1}$ for $i<n$. Or, in the other notation:
\[
    J(x_1,\cdots,x_n)=[x_1,\cdots,x_n].
\]

\begin{defn}\label{Z:def:Z+}
Let $\Z_+$ denote the category with both linearly and cyclically
ordered sets (ie, $Ob(\Z_+)=Ob(\nL)\coprod Ob(\Z)$ and three
kinds of morphisms:
\begin{enumerate}
  \item the usual morphisms (degree 1 monomorphisms) between
  objects of $\Z$,
  \item the usual morphisms (order preserving monomorphisms) between
  objects of $\nL$ and
  \item a morphism $f\co B\to C$ from a linearly ordered set $B$ to
  a cyclically ordered set $C$ is defined to be a
  morphism $f\co J(B)\to C$ in $\Z$ (and composition is given
  by $f\circ g=fJ(g)\co J(A)\to C$ if $g\co A\to B$ is a morphism in
  $\nL$).
\end{enumerate}
There are no morphisms from $\Z$ to $\nL$.
\end{defn}

Note that $\nL, \Z$ are full subcategories of $\Z_+$. A retraction
$\overline{J}\co \Z_+\to \Z$ is given by $J$ on $\nL$ and the
identity on $\Z$. The identity map on sets $C\to \overline{J}(C)$
is a natural transformation from the identity functor on $\Z_+$ to
the functor $\overline{J}$. This proves the following.

\begin{prop}\label{Z:prop:Z simeq Z+}
$\Z$ is a deformation retract of $\Z_+$.
\end{prop}

The difference between $\Z$ and $\Z_+$ is that $\Z_+$ has a base
point up to homotopy. The full subcategory of $\nL$ in $\Z_+$ is
contractible and therefore serves as a homotopy base point for the
category $\Z_+$. (A contraction is given by adding one point on
the left then delete all the other points.)

\subsection{Homotopy type of $|\Z|\simeq|\Z_+|$}

Although $\Z_+$ may be unfamiliar, the homotopy type of $\Z$ is
well-known. See, eg, \cite{[IK1:Borel2]} or \cite{[I:Mult-HC]}.

\begin{thm}
$|\Z|\simeq|\Z_+|\simeq\CC P^\infty$.
\end{thm}

The universal circle bundle over $\CC P^\infty$ pulls back to a
circle bundle over the geometric realization $|\Z_+|$ of $\Z_+$
given by $S^1(C,\s)$ over $(C,\s)$ with or without base point. A
precise construction will be given later.

\subsection{Powers of the Euler class}\label{Z:subsec:powers of
eZ}

The Euler class of the universal circle bundle over $\Z$ is a $2$
dimensional integral cohomology class $e_\Z\in H^2(\Z,\ZZ)$. It is
represented by a $2$--cocycle which assigns an integer to every
$2$--simplex $A\to B\to C$ in $\Z$. However, the classical method
is to choose a connection and integrate its curvature. This
procedure, carried out in the next section, produces a rational
$2$--cocycle on $\Z$ which represents this integral class. Since
$\CC P^\infty$ has no torsion in its homology, the integral class
is uniquely determined by this rational cocycle.

The cocycle representing the Euler class can be described as
follows. Given a $2$--simplex $A\to B\to C$ in $\Z$, we choose
elements $a,b,c$ in $A,B,C$. The \emph{sign} of $(a,b,c)$ is $+$
if the images of $a,b,c$ in $C$ are distinct and in cyclic order.
The sign is $-1$ if they are distinct and in reverse cyclic order.
If they are not distinct then the sign is $0$. The cocycle
$c_\Z(A,B,C)$ is defined to be $-\frac12$ times the expected value
of this sign.

For example, if $A=[a],B=[a,b],C=[a,c,b]$ with morphisms being
inclusion maps then the probability is $\frac16$ that distinct
elements of $A,B,C$ will be chosen. The sign of $(a,b,c)$ is
negative since it is an odd permutation of the given cyclic
ordering of $C$ so
\[
    c_\Z(A,B,C)=-\left(-\frac12\right)\left(\frac16\right)=+\frac1{12}.
\]
More generally, suppose that
\[
  C_\ast=(C_0\to C_1\to\cdots\to C_{2k})
\]
is a $2k$ simplex in $\Z$. Then the \emph{cyclic set cocycle}
$c_\Z^k$ is defined on $C_\ast$ by
\begin{equation}\label{Z:eq:formula for cZk}
  c_\Z^k(C_\ast)=
  \frac{(-1)^k k!\sum \sgn(a_0,a_1,\cdots,a_{2k})}
  {(2k)!|C_0|\cdots|C_{2k}|}
\end{equation}
where the sum is taken over all $a_i$ in the image
of $C_i$ in $C_{2k}$ for $i=0,\cdots,2k$ and the sign of
$(a_0,a_1,\cdots,a_{2k})$ is given by comparing this ordering with
the ordering induced by the cyclic ordering of $C_{2k}$. (The sign
is zero if these elements are not distinct.)

Note that in (\ref{Z:eq:formula for cZk}) the sum is the same if
we take only those $2k+1$ tuples $(a_0,a_1,\cdots,a_{2k})$ where
$a_i$ is in the image of $C_i$ but not in the image of $C_{i-1}$.
(Otherwise, take $i$ minimal so that $a_i\in C_j$ for some $j<i$
and switch $a_i$ and $a_j$, where $j$ is minimal. This gives
another summand with the opposite sign and the described operation
is an involution on the set of summands that we are deleting.)

\begin{prop}\label{Z:prop:cZk is a cocycle}
The cyclic set cocycle $c_\Z^k$ is a rational cocycle on $\Z$.
\end{prop}

\begin{proof}
Suppose that $C_0\to \cdots\to C_{2k+1}$ is a $2k+1$ simplex in
$\Z$. Choose one element from each $C_i$ at random with equal
probability. Let $a_i$ be the image of this element in $C_{2k+1}$.
Then the following alternating sum of the signs vanishes.
\[
    \sum_{i=0}^{2k+1}(-1)^i
    \sgn(a_0,\cdots,\widehat{a_i},\cdots,a_{2k+1})=0
\]
This is obvious if these elements are in cyclic order in
$C_{2k+1}$. It is also clear that if two consecutive elements are
reversed in the cyclic ordering then every summand in the above
expression changes sign. Finally, if the elements $a_i$ are not
distinct then all summands are zero except for two terms of
opposite sign.

Since expected value is a linear function, the sum of the expected
values of the summands is zero, ie, the expected value of the
sign is a cocycle.
\end{proof}

\subsection{Extension of $c_\Z^k$ to $\Z_+$}

For the purpose of constructing the combinatorial
Miller--Morita--Mumford\break classes $[c_\Fat^k]$ we only need $c_\Z^k$
on cyclically ordered sets. However, when we go to framed graphs
we will need to extend this to $\Z_+$ in such a way that it is
identically zero on $\nL$.

First we note that the pull-back along the functor $J\co \nL\to\Z$ of
the cocycle $c_\Z^k$ must be a coboundary since $\nL$ is
contractible. It is in fact the coboundary of the $2k-1$ cochain
$s_k$ given on $C_\ast=(C_0\to\cdots\to C_{2k-1})$ by
\begin{equation}\label{Z:eq:the 2k-1 cochain sk}
  s_k(C_\ast)=
  \frac{(-1)^kk!\sum \sgn(a_0,a_1,\cdots,a_{2k-1})}
  {(2k)!|C_0|\cdots|C_{2k-1}|}.
\end{equation}
As before, we note that this sum is unchanged if we delete terms
where $a_i\in C_{i-1}$ for some $i$.

\begin{defn}\label{Z:def:czk on Z+}
Given a $2k$--simplex $C_\ast=(C_0\to\cdots\to C_{2k})$ in $\Z_+$
we define $c_{\Z^+}^k(C_\ast)$ as follows.
\begin{enumerate}
  \item[(a)] If the last two objects lie in $\Z$ we let
  $c_{\Z^+}^k(C_\ast)=c_\Z^k\overline{J}(C_\ast)$ where
  $\overline{J}\co \Z_+\to\Z$ is the retraction which is equal to $J$
  on $\nL$.
  \item[(b)] If all the objects lie in $\nL$ then we let
  $c_{\Z^+}^k(C_\ast)=0$.
  \item[(c)] If $C_{2k-1}\in\nL$ and $C_{2k}\in\Z$ then
  let
\begin{equation}\label{Z:eq:formula for cZk on Z+}
    c_{\Z^+}^k(C_\ast)=c_\Z^k\overline{J}(C_\ast)-s_k(C_0\to\cdots\to
    C_{2k-1})
\end{equation}
where $s_k$ is given in (\ref{Z:eq:the 2k-1 cochain sk}) above.
\end{enumerate}
\end{defn}
Note that in all three cases $c_{\Z^+}^k(C_\ast)=0$ if the sets
$C_0,\cdots,C_{2k-1}$ do not have distinct cardinalities.

Since $\delta s_k=c_\Z|\nL$ it follows that:

\begin{prop}\label{Z:prop:cZk is a cocycle on Z+}
The extended cyclic set cocycle $c_{\Z^+}^k$ is a rational
$2k$--cocycle on $\Z_+$.
\end{prop}

In the next section we prove the following.

\begin{thm}\label{Z:thm:cZk represents eZk on Z}
The rational $2k$--cocycle $c_\Z^k$ represents the $k^{\rm th}$ power
$e_\Z^k$ of the Euler class $e_\Z\in H^2(\Z,\ZZ)$.
\end{thm}

Since $|\Z|\simeq |\Z^+|$ and $c_\Z^k=c_{\Z^+}^k|\Z$ we get:

\begin{cor}\label{Z:cor:cZk+ represents eZk+ on Z+}
The extended cyclic set cocycle $c_{\Z^+}^k$ represents the $k^{\rm th}$
power $e_{\Z^+}^k$ of the Euler class $e_{\Z^+}$ of ${\Z^+}$.
\end{cor}

%
%

\section{The curvature form on $|\Z|$}  

In this section we obtain the cyclic $2k$--cocycle $c_\Z^k$ as an
elementary exercise in differential geometry. Briefly the idea is
that we want to find a natural connection $A$ on the geometric
realization $|\Z|$ of the category $\Z$, take the powers of the
associated curvature form
\[
    \Omega=dA
\]
and integrate over the even dimensional simplices of $|\Z|$. The
same curvature form appears in \cite{[Kontsevich:Airy]} giving the
Euler class for a space $BU(1)^{comb}$ which is closely related to
$|\Z|$.

\subsection{Smooth families of cyclically ordered sets}

The first step is to construct a (piecewise) smooth space of
cyclically ordered sets. The idea is simple. We view a cyclically
ordered set with $n$ elements as being $n$ points evenly spaced on
a unit circle. A smooth version of this is to divide the circle
into $n$ arcs of varying length (but with constant total length
$2\pi$). By letting the lengths of some arcs go to zero we can
change the number of elements in the cyclically ordered set in a
continuous way.

We define the \emph{weight} of an arc to be its length divided by
$2\pi$. Then, up to rotation, an element in this space is
represented by a cyclically ordered sequence of $n$ nonnegative
real numbers $[w_1,\cdots,w_n]$ whose sum is $1$. We make this
precise:

\begin{defn}\label{dZ:def:cyclic weight system}
A \emph{cyclic weight system} is a triple $(C,\s,w)$ where
$(C,\s)$ is a cyclically ordered set and $w\co C\to I$ is a
nonnegative real valued function on $C$ so that \[\sum_{x\in
C}w(x)=1.\]
\end{defn}

Given a cyclic weight system $(C,\s,w)$, the \emph{canonical
circle} over $(C,\s,w)$ is given by
\[
    S^1(C,\s,w)=\coprod_{x\in C}x\times[0,w(x)]/\sim
\]
where the identifications are given by $(x,w(x))\sim(\s(x),0)$,
ie, the line segments $x\times[0,w(x)]$ are connected end to end
in a circle. If $w_t$, $t\in P$, is a smooth family of weights on
a fixed cyclically ordered set $(C,\s)$ then we can form a smooth
principal $S^1$--bundle over $P$ by:
\[
    S^1(C,\s,w_P)=\coprod_{x\in C}\{(s,t)\in I\times P\st s\leq
    w_t(x)\}/\sim
\]
with the fiberwise identification $(x,w_t(x),t)\sim(\s(x),0,t)$.

If we choose a starting point, the elements of $C$ can be written
$(x_1,\cdots,x_n)$ (with cyclic order $[x_1,\cdots,x_n]$) and we
can write $w_j=w(x_j)$. If we place line segments of length $w_j$
next to each other in sequence on the real line starting at the
origin then the center of mass of the $j^{\rm th}$ segment will be
located at a point
\[
    s_j=w_1+\cdots+\w_{j-1}+\frac12 w_j
\]
units from the starting point. On the circle this will be the
point $\exp(2\pi is_j)$.

The \emph{angular momentum} of the circle is then given by:
\[
    \sum w_j\exp(-2\pi is_j)d\exp(2\pi is_j)=2\pi i\sum w_jds_j.
\]
This means that, with respect to the inertial frame, our
coordinate system is rotating in the other direction at this rate.
Consequently the \emph{inertial connection} (in this coordinate
frame) is given by
\[
    A=-2\pi i\sum w_jds_j=-2\pi i\sum_{1\leq i<j\leq n}w_jdw_i
    -\pi i \sum w_jdw_j.
\]
The \emph{curvature} of this form is
\[
    \Omega=dA= dA=2\pi i\sum_{1\leq i<j\leq n}dw_i\wedge
    dw_j.
\]
Therefore, the \emph{Euler class} of the canonical circle bundle
over the space of cyclic weight systems is given by the
differential form
\begin{equation}\label{dZ:eq:Euler form}
  \frac{i}{2\pi}\Omega=-\sum_{1\leq i<j\leq n}dw_i\wedge
    dw_j.
\end{equation}
This is invariant under cyclic permutation of the $w_j$ since the
summands with $i=1$ add up to
$dw_1\wedge(dw_2+\cdots+dw_n)=dw_1\wedge(-dw_1)=0$. Similarly, the
terms with $j=n$ add up to zero.

We interpret the ``space of cyclic weight systems'' to be the
geometric realization $|\Z|$ of $\Z$.

In \cite{[Kontsevich:Airy]} Konsevich obtained (\ref{dZ:eq:Euler
form}) (with all terms having $j=n$ deleted) as the Euler class on
a space $BU(1)^{comb}$ which he defined to be the space of
isomorphism classes of cyclic weight systems. If $\Z_0$ is the
full subcategory of $\Z$ given by choosing one object from each
isomorphism class then we get a quotient map
$|\Z_0|\twoheadrightarrow BU(1)^{comb}$ which is a rational
homotopy equivalence. Kontsevich shows that the Euler class of
$BU(1)^{comb}$ is given by a $2$--form $\omega$. The differential
form (\ref{dZ:eq:Euler form}) is the pull-back of this $2$--form.

\subsection{Simplices in $|\Z|$}

Strictly speaking $|\Z|$ is the geometric realization of the
\emph{simplicial nerve} $\N\bu(\Z)$ of $\Z$. A $p$--simplex
\[
    C_\ast=(C_0\to C_1\to\cdots\to C_p)
\]
in $\Z$ is one element of $\N_p(\Z)$ but it represents a geometric
$p$--simplex
\[
    \Delta^p\times\{C_\ast\}\subseteq |\Z|.
\]
By definition, $|\Z|$ is the union of these geometric simplices:
\[
    |\Z|=\coprod_{x\in\N_p(\Z),p\geq0}\Delta^p\times
    x/\sim.
\]
The vertices $v_0,\cdots,v_p$ of $\Delta^p$ correspond to the
cyclic sets $C_0,\cdots,C_p$. The other points $t\in\Delta^p$
correspond to cyclic weight systems given by mass functions on
$C_p$. The set $C_j$ is identified with the mass function on $C_p$
which is $1$ on the image of $C_j$ and $0$ in the complement. For
simplicity of notation we will identify $C_j$ with its image in
$C_p$.

We parametrize the $p$--simplex $\Delta^p$ by
\[
    \Delta^p=\{t\in\RR^p\st 1\geq t_1\geq t_2\geq\cdots\geq
    t_p\geq0\}
\]
and we take $t_0=1,t_{p+1}=0$ to be fixed. Then the cyclic weight
system $C(t)$ for $t\in \Delta^p$ will be given by the \emph{mass
function}
\[\mu_t\co C_p\to I\]
given by $\mu_t(a)=t_j$ if $a\in C_j-C_{j-1}$. The weights are the
normalized masses
\[
    w_t(a)=\frac{\mu_t(a)}{M_t}=\frac{\mu_t(a)}{\sum\mu_t(a)}
\]
where $M_t=\sum\mu_t(a)$ is the total mass. The weights are
ordered according to the cyclic ordering of the elements of $C_p$.

The \emph{canonical circle bundle} $E_\Z$ over $|\Z|$ is given by
\[
    E_\Z=\coprod_{x\in\N_p(\Z),p\geq0}S^1(C_p,w_{\Delta^p})\times x/\sim.
\]
Note that each piece $S^1(C_p,w_{\Delta^p})\times x$ is a smooth
principal $S^1$ bundle over $\Delta^p\times x$. Consequently,
$E_\Z$ is a piecewise smooth principal $S^1$--bundle over $|\Z|$.

The $j^{\rm th}$ vertex $v_j$ of $\Delta^p$ is given by
\[
    t_0=t_1=\cdots=t_j=1,\quad t_{j+1}=\cdots=t_p=t_{p+1}=0.
\]
This agrees with the discussion above since it assigns a mass of
$1$ to the elements of $C_j$. The barycentric coordinates on
$\Delta^p$ are given by $t_j'=t_j-t_{j+1}$ so that the $j^{\rm th}$ face
is given by $t_j=t_{j+1}$.

\subsection{The Euler class on $2$--simplices}

Now take $p=2$. Take a fixed $2$--simplex $C_0\to C_1\to C_2$. Let
$a=|C_0|,b=|C_1|-a$ and $c=|C_2|-a-b$. Denote the elements of
$C_2$ in cyclic order by $(x_1,x_2,\cdots,x_n)$ where $n=a+b+c$.

If $t=(t_1,t_2)\in\Delta^2$ then the mass function $\mu_t$ is
given by
\[
    \mu_t(x_i)=
  \begin{cases}
    t_0=1 & \text{if $x_i\in C_0$}, \\
    t_1 & \text{if $x_i\in C_1-C_0$}, \\
    t_2 & \text{if $x_i\in C_2-C_1$}.
  \end{cases}
\]
We note that $t_0$ occurs $a$ times, $t_1$ occurs $b$ times and
$t_2$ occurs $c$ times. Thus the total mass is
\[
    M_t=\sum \mu_t(x_i)=a+bt_1+ct_2.
\]
The weight (relative mass) of $x_i$ is
\[
    w_t(x_i)=\frac{\mu_t(x_i)}{M_t}=\frac{t_j}
    {
        a+bt_1+ct_2
    }
\]
where $j=0,1$ or $2$.

\begin{thm}\label{dZ:thm:the Euler class}
The $2$--cocycle on $\Z$ whose value on the $2$--simplex $(C_0\to
C_1\to C_2)$ is given by
\[
    \int_{\Delta^2}\frac{i}{2\pi}\Omega=
    -\int_{1\geq t_1\geq t_2\geq0}
    \sum_{1\leq i< j\leq n}
    dw_t(x_i)\wedge
    dw_t(x_j)
\]
is equal to the $2$--cocyle $c_\Z$ of section~\ref{Z:subsec:powers
of eZ}.
\end{thm}

\proof
Up to sign there are only three possibilities for the $2$--form $
dw_t(x_i)\wedge dw_t(x_j)$:
\begin{align}\label{dZ:eq:term1 in iOmega}
    d\left(
        \frac1{M_t}
        \right)
    \wedge
    d\left(
        \frac{t_1}{M_t}
        \right)
&=
    -\frac{dM_t}{M_t^2}\wedge\frac{dt_1}{M_t}
    =
    \frac{-cdt_2\wedge dt_1}{M_t^3}
    =
    c\frac{dt_1\wedge dt_2}{M_t^3}
\\
\label{dZ:eq:term2 in iOmega}
    d\left(
        \frac1{M_t}
        \right)
    \wedge
    d\left(
        \frac{t_2}{M_t}
        \right)
&=
    -\frac{dM_t}{M_t^2}\wedge\frac{dt_2}{M_t}
    =
    -b\frac{dt_1\wedge dt_2}{M_t^3}
\\
\label{dZ:eq:term3 in iOmega}
    d\left(
        \frac{t_1}{M_t}
        \right)
    \wedge
    d\left(
        \frac{t_2}{M_t}
        \right)
&=
    \left[
        1-\frac{bt_1}{M_t}-\frac{ct_2}{M_t}
    \right]
    \frac{dt_1\wedge dt_2}{M_t^2}
    =
    a    \frac{dt_1\wedge dt_2}{M_t^3}
\end{align}
We interpret (\ref{dZ:eq:term1 in iOmega}) as a sum of $c$ terms
(one for each $x_k\in C_2-C_1$) and similarly for
(\ref{dZ:eq:term2 in iOmega}) and (\ref{dZ:eq:term3 in iOmega}).
Then for every triple of indices $(i,j,k)$ so that $x_i\in C_0$,
$x_j\in C_1-C_0$ and $x_k\in C_2-C_1$ we get three terms, one of
each kind, adding up to:
\[
\left( \sgn(j-i)-\sgn(k-i)+\sgn(k-j)
    \right)
    \frac{dt_1\wedge dt_2}{M_t^3}.
\]
A little thought will show that the sum of signs is $1$ if $i,j,k$
are in cyclic order and $-1$ if not. Ie,
\begin{equation}\label{dZ:eq:sign sum with three terms}
    \sgn(j-i)-\sgn(k-i)+\sgn(k-j)=\sgn(i,j,k).
\end{equation}
Furthermore, we have the easy double integral:
\begin{equation}\label{dZ:eq:double integral}
    \int_{1\geq t_1\geq t_2\geq0}\frac{dt_1dt_2}
    {(a+bt_1+ct_2)^3}
    =
    \frac1{2a(a+b)(a+b+c)}.
\end{equation}
Putting these together we get:
$$
    \int_{\Delta^2}\frac{i}{2\pi}\Omega
    =
   - \sum_{i,j,k\in C_2}\frac{
        \sgn(i,j,k)
    }{
        2a(a+b)(a+b+c)
    }
    =
    c_\Z(C_0\to C_1\to C_2)\eqno{\qed}
$$

\subsection{$e_\Z^k$ on $2k$--simplices}

Now take $p=2k$. Let $C_0\to\cdots\to C_{2k}$ be a $2k$--simplex in
$\Z$. Let $a_0=|C_0|$ and $a_j=|C_j|-|C_{j-1}|$ for $j\geq1$. Let
$(x_1,\cdots,x_n)$ denote the elements of $C_{2k}$ in cyclic
order.

The mass function $\mu_t$ for $t\in \Delta^{2k}$ is given by
$\mu_t(x_i)=t_j$ if $x_i\in C_j-C_{j-1}$. Then the total mass is
\[
    M_t=\sum_{i=1}^n\mu_t(x_i)=\sum_{j=0}^{2k}a_jt_j.
\]
The weight of $x_i$ is
\[
    w_t(x_i)=\frac{\mu_t(x_i)}{M_t}=\frac{t_j}{a_0t_0+\cdots+a_{2k}t_{2k}}
\]
for some $j$.

 A cocycle for $e_\Z^k$ is given on $2k$--simplices by
integrating the $2k$--form $\left(\frac{i}{2\pi}\Omega\right)^k$
over $\Delta^{2k}$. This form can be expanded:
\begin{align*}
  \left(\frac{i}{2\pi}\Omega\right)^k
     &=
    (-1)^k\sum_{i_1<j_1}\cdots\sum_{i_k<j_k}dw_t(x_{i_1})\wedge
    dw_t(x_{j_1})\wedge\cdots\wedge dw_t(x_{i_k})\wedge
    dw_t(x_{j_k})
\\
    &=_{(a)}(-1)^kk!\sum_{i_1<i_2<\cdots<i_k}\sum_{i_1<j_1}\cdots\sum_{i_k<j_k}
    dw_t(x_{i_1})\wedge \cdots\wedge
    dw_t(x_{j_k})
\\
    &=_{(b)}(-1)^kk!\sum_{i_1<j_1<i_2<\cdots<i_k<j_k}
    dw_t(x_{i_1})\wedge \cdots\wedge
    dw_t(x_{j_k})
\\
    &=(-1)^kk!\sum_{i_1<i_2<\cdots<i_{2k}}
    dw_t(x_{i_1})\wedge \cdots\wedge
    dw_t(x_{i_{2k}})
\end{align*}
where (a) is by symmetry and (b) follows from the fact that the
summands in the second line not in the third come in cancelling
pairs: Take the first $j_p>i_{p+1}$ and switch it with $j_{p+1}$.

Suppose that $w_t(x_{i_p})={t_{j_p}}/{M_t}$. Then there are
basically only two possibilities for $dw_t(x_{i_1})\wedge
\cdots\wedge dw_t(x_{i_{2k}})$:
\begin{enumerate}
  \item[(a)] If $j_p\noteq0$ for all $p$ then
  \[
    dw_t(x_{i_1})\wedge \cdots\wedge
    dw_t(x_{i_{2k}})
    =
    \sgn(j_1,\cdots,j_{2k})
    d\left(\frac{t_1}{M_t}\right)\wedge\cdots\wedge
    d\left(\frac{t_{2k}}{M_t}\right)
  \]
  \[
    =\sgn(j_1,\cdots,j_{2k})
    \left(
    1-\frac{a_1t_1}{M_t}-\cdots-\frac{a_{2k}t_{2k}}{M_t}
    \right)
    \frac{dt_1\wedge dt_2\wedge\cdots\wedge dt_{2k}}
    {M_t^{2k}}
  \]
  \[
    =a_0\sgn(0,j_1,\cdots,j_{2k})
    \frac{dt_1\wedge dt_2\wedge\cdots\wedge dt_{2k}}
    {M_t^{2k+1}}
  \]
  \item[(b)] If $j_p=0$ for some $p$ then the corresponding weight is
   $w_t(x_{i_p})=\frac1{M_t}$ and we have
   $d\left(\frac1{M_t}\right)$ instead of
   $d\left(\frac{t_q}{M_t}\right)$ for some $q$. This means we
   should replace the term $d\left(\frac1{M_t}\right)$ by
   $-\frac{a_qdt_q}{M_t^2}$ and we get:
   \[
    dw_t(x_{i_1})\wedge \cdots\wedge
    dw_t(x_{i_{2k}})
    =
    a_q\sgn(q,j_1,\cdots,j_{2k})
    \frac{dt_1\wedge \cdots\wedge dt_{2k}}
    {M_t^{2k+1}}
   \]
\end{enumerate}

As before, we interpret (a) as a sum of $a_0$ terms, one for each
element of $C_0$, and similarly for (b). Then for every choice of
$x_{i_0}\in C_0, x_{i_1}\in C_1-C_0,\cdots, x_{i_{2k}}\in
C_{2k}-C_{2k-1}$ we get the following.
\begin{equation}\label{dZ:eq:terms coming from x0 to x2k}
  \sum_{q=0}^{2k}\sgn(x_{i_q},x_{i_0},\cdots,\widehat{x_{i_q}},\cdots,x_{i_{2k}})
    \frac{dt_1\wedge \cdots\wedge dt_{2k}}
    {M_t^{2k+1}}
\end{equation}
The permutations in (\ref{dZ:eq:terms coming from x0 to x2k}) are
the inverses of the permutations in (a),(b) above so they have the
same sign.

The summands in (\ref{dZ:eq:terms coming from x0 to x2k}) are
equal with alternating signs. Since there are an odd number of
terms it is equal to its first summand. The $2k$--form
$(1/M_t^{2k+1})dt_1\wedge\cdots\wedge dt_{2k}$ has integral:
\begin{align}\label{dZ:eq:multiple integral line1}
  \int_{\Delta^{2k}}
  \frac{dt_1\wedge \cdots\wedge dt_{2k}}
    {M_t^{2k+1}}
&=
    \int_{1\geq t_1\geq\cdots\geq t_{2k}\geq0}
    \frac{dt_1\cdots dt_{2k}}
    {(a_0+a_1t_1+\cdots+a_{2k}t_{2k})^{2k+1}}\\
\label{dZ:eq:multiple integral line2}
&=
    \frac1{(2k)!a_0(a_0+a_1)\cdots(a_0+\cdots+a_{2k})}
\end{align}
which is an easy induction on $2k$.

Combining this with the formula for
$\left(\frac{i}{2\pi}\Omega\right)^k$ we get:

\begin{thm}\label{dZ:thm:equation for cZk}
The integral of $\left(\frac{i}{2\pi}\Omega\right)^k$ over the
$2k$--simplex $\Delta^{2k}\times\{C_\ast\}$ is
\begin{equation}\label{dZ:eq:cZk as an integral}
  \int_{\Delta^{2k}}\left(\frac{i}{2\pi}\Omega\right)^k
  =
  (-1)^kk!\frac{\sum \sgn(x_{i_0},\cdots,x_{i_{2k}})}
    {(2k)!|C_0|\cdots|C_{2k}|}
\end{equation}
where the sum is taken over all $x_{i_0}\in C_0,x_{i_1}\in
C_1-C_0,\cdots,x_{i_{2k}}\in C_{2k}-C_{2k-1}$.
\end{thm}

In other words, the deRham cocycles for the powers of the Euler
class on $|\Z|$ are equal to the combinatorial cocycles $c_\Z^k$.

%
%

\section{Combinatorial formula for MMM classes}

We construct cocycles $c_\Fat^k$ on the category of fat graphs
$\Fat$ by evaluating the cyclic set cocycles $c_\Z^k$ on each
vertex. The classifying space of this category is well-known to be
homotopy equivalent to the disjoint union of classifying spaces of
mapping class groups $M_g^s$ of surfaces of genus $g$ with $s\geq
max(1,3-2g)$ punctures.
\[
    |\Fat|\simeq\coprod_{s\geq max(1,3-2g)}BM_g^s
\]
(This is Theorem \ref{FG:thm:FGgs is BMgs} below.) Thus a
cohomology class for $\Fat$ gives a cohomology class for each
mapping class group $M_g^s$.

By direct computation we show that the cohomology classes
$[c_\Fat^k]$ are dual to the Witten cycles $W_{k}$. More
precisely,
\[
    [c_\Fat^k]=(-1)^k\frac{k!}{(2k+1)!}[W_{k}]^\ast.
\]
In Theorem \ref{Calc:thm:proportionality constant between MMM
and b} below we will show that $[c_\Fat^k]=-2\wt{\k}_k$. So the
\emph{adjusted fat graph cocycle} $-\frac12c_\Fat^k$ is a
combinatorial formula for $\wt{\k}_k$.

\subsection{The category of fat graphs}

We define a \emph{fat graph} to be a finite connected graph
$\Gamma$ possibly with loops and multiple edges in which every
vertex has valence $\geq3$ together with a cyclic ordering of the
edges incident to each vertex. To be precise and to fix our
notation, the fat graph $\Gamma$ consists of
\begin{enumerate}
  \item $\Gamma^0$, the set of \emph{vertices},
  \item $\Gamma^\frac12$, the set of \emph{half edges},
  \item $\d\co \Gamma^\frac12\to\Gamma^0$, the \emph{incidence} or
  \emph{boundary} map so that $|\d^{-1}(v)|\geq3$ for all vertices
  $v$,
  \item $\s$, a cyclic ordering on each set $\d^{-1}(v)$ and
  \item $a\mapsto\overline{a}$, a fixed-point free involution
  on $\Gamma^\frac12$ whose orbits we call \emph{edges}.
\end{enumerate}
Note that each edge $\{a,\overline{a}\}$ has two orientations
$e=(a,\overline{a})$ and $\overline{e}=(\overline{a},a)$. Each
oriented edge $e=(a,b)$ has a \emph{source} $s(e)=\d a$ and
\emph{target} $t(e)=\d b$.

For several reasons we need to consider the cyclically ordered set
of ``angles'' between incident half edges at each vertex of a fat
graph. An \emph{angel} (at $v$) is defined to be an ordered pair
of half edges $(a,b)$ so that $\d a=\d b=v$ and $b=\s(a)$. In
other words, $a,b$ are incident to the same vertex $v$ and $b$ is
one step counterclockwise from $a$. Let $C(v)$ be the set of
angles at $v$. Then $C(v)$ has a cyclic ordering
$\s(a,b)=(b,\s(b))$.

For example, the figure ``$\infty$'' has one vertex $v$ of valence
$4$ with $\d^{-1}(v)=[a,\overline{a},b,\overline{b}]$ and
$C(v)=[(a,\overline{a}),(\overline{a},b),
    (b,\overline{b}),(\overline{b},a)]$. As before, we denote cyclically ordered sets by square brackets.

A fat graph $\Gamma$ is evidently equal to the core of some connected
oriented punctured surface $\Sig_\Gamma$ which is well-defined up
to homeomorphism.

A \emph{morphism} $f\co \Gamma_1\to\Gamma_2$ of fat graphs is a
morphism of graphs where the inverse image of every open edge is
an open edge and the inverse image of every vertex is a tree in
$\Gamma_1$ with the cyclic ordering of the half edges incident to each
vertex of $\Gamma_2$ corresponds to the cyclic ordering of the half
edges incident to the corresponding tree in $\Gamma_1$. (In other
words, the surfaces are homeomorphic.)

The \emph{codimension} of a graph is defined to be the non-negative integer
\[
    \codim\Gam=\sum( \val(v)-3).
\]
Thus, $\codim\Gam=0$ if and only if $\Gam$ is trivalent. It is important to note that, for any morphism $f\co \Gam_1\to\Gam_2$ which is not an isomorphism, $\codim \Gam_1>\codim\Gam_2$.

The category of all fat graphs will be denoted $\Fat$. Since the
punctured surface $\Sig_\Gamma$ is fixed up to homeomorphism on
each component of $\Fat$ we have:
\[
    \Fat=\coprod \Fat_g^s
\]
where $\Fat_g^s$ is the full subcategory of fat graphs $\Gamma$ so
that $\Sig_\Gamma$ is a surface of genus $g$ with $s$ punctures.

There is a well-known correspondence between fat graphs, the
moduli space of curves and the mapping class group. In the present
context it says the following.

\begin{thm}[Penner \cite{[Penner87:decorated]}, Strebel \cite{[Strebel]}]\label{FG:thm:FGgs is BMgs}
The geometric realization of the cate\-gory $\Fat_g^s$ is homotopy
equivalent to the classifying space of the mapping class group
$M_g^s$ of genus $g$ surfaces with $s$ punctures provided that
$\Fat_g^s$ is non-empty, ie, $s\geq1$ and $s+2g\geq3$.
\end{thm}

\begin{proof} (For a more detailed proof see \cite{[I:BookOne]},
Theorem 8.6.3.) Let $S$ be a fixed oriented surface of genus $g$ with
$s$ boundary components. Then, by a theorem of Culler and Vogtmann
\cite{[Culler-Vogtmann-86]}, the space of all pairs $(\Gamma,f)$
where $\Gamma$ is a fat graph (an element of $|\Fat|$) and $f$ is
an orientation preserving homeomorphism $f\co \Sig_\Gamma\to S$ is
contractible and $\Homeo_+(S)$ acts freely on this space with
quotient $|\Fat_g^s|$. Thus $|\Fat_g^s|\simeq B\Homeo_+(S)\simeq
BM_g^s.
$
\end{proof}

\subsection{The fat graph cocycle $c_{\Fat}^k$}

For each vertex $v$ of $\Gamma_1$, a morphism
$f\co \Gamma_1\to\Gamma_2$ sends the angle set $C(v)$ monomorphically
into $C(f(v))$ in a cyclic order preserving way. Thus to a
$2k$--simplex
\[
    \Gamma_\ast=(\Gamma_0\to\Gamma_1\to\cdots\to\Gamma_{2k})
\]
in the nerve of $\Fat$ we can extract several $2k$--simplices in
the nerve of $\Z$, one for each vertex of $\Gamma_0$.

\begin{defn}\label{FG:def:the combinatorial cFGk}
Let $c_\Fat^k$ be the $2k$--cochain on $\Fat$ given by
\[
    c_\Fat^k(\Gamma_\ast)=\sum_{v\in\Gamma_0^0} m(v)c_\Z^k(C(v)\to
    C(f_1(v))\to\cdots\to C(f_{2k}(v)))
\]
where $m(v)=\val(v)-2$ is the \emph{multiplicity} of $v$ and
$f_i=f_{i0}\co \Gamma_0\to\Gamma_i$ is the composition
\[
    f_{i0}\co \Gamma_0\xrarrow{f_{10}}\Gamma_1\xrarrow{f_{21}}\cdots
    \xrarrow{f_{i-1\ i}}\Gamma_i
\]
of arrows in $\Gamma_\ast$.
\end{defn}

Every time an edge collapses, two half edges disappear.
Consequently, the multiplicity of the resulting vertex is the sum
of the multiplicities of the original two vertices. More
generally, we have:

\begin{lem}\label{FG:lem:multiplicities add up}
Given any morphism $f\co \Gamma_1\to\Gamma_2$ in $\Fat$ and any
vertex $v$ in $\Gamma_2$ we have:
\[
    m(v)=\sum_{w\in f^{-1}(v)}m(w),
\]
ie, the multiplicity of $v$ is the sum of the multiplicities of
the vertices which collapsed to $v$.
\end{lem}

\proof
This follows from the fact that $T=f^{-1}(v)$ is a tree. Thus $T$
has $n$ edges and $n+1$ vertices $w_0,\cdots,w_n$. So,
$$
    m(v)=\val(v)-2=\sum\val(w_i)-2n-2=\sum m(w_i).\eqno{\qed}
$$

\begin{thm}\label{FG:thm:cFGk is a cocycle for an integral class}
$c_\Fat^k$ is a rational $2k$--cocycle on $\Fat$ which determines a
well-defined integral cohomology class
\[
    [c_\Fat^k]\in H^{2k}(\Fat,\ZZ).
\]
\end{thm}

\begin{proof}
Given any $2k+1$--simplex
$\Gamma_\ast=(\Gamma_0,\cdots,\Gamma_{2k+1})$ we have:
\begin{align*}
  \delta c_\Fat^k(\Gamma_\ast) & =
    \sum_{i=0}^{2k+1}(-1)^i
    c_\Fat^k(\Gamma_0,\cdots,\widehat{\Gamma_i},\cdots,\Gamma_{2k+1})
    \\
    &
    =\sum_{v_1\in\Gamma_1^0}m(v_1)c_\Z^k(C(v_1),\cdots,C(v_{2k+1}))\\
&+
    \sum_{i=1}^{2k+1}(-1)^i\sum_{v_0\in\Gamma_0^0}m(v_0)
    c_\Z^k(C(v_0),\cdots,\widehat{C(v_i)},\cdots,C(v_{2k+1}))
\end{align*}
where the $v_i$ are related by $v_i=f_{ij}(v_j)$ for all $j<i$.
Since $c_\Z^k$ is a cocycle, the last sum is equal to
\[
    -\sum_{v_0\in\Gamma_0^0}m(v_0)
    c_\Z^k(C(v_1),\cdots,C(v_{2k+1}))
\]
which exactly cancels the second sum by
Lemma~\ref{FG:lem:multiplicities add up}. Thus $c_\Fat^k$ is a
(rational) cocycle. But the above argument uses only the fact that
$c_\Z^k$ is a $2k$--cycle on the category $\Z$. Therefore we may
replace $c_\Z^k$ with an integral cocycle. Since $|\Z|\simeq\CC
P^\infty$ has no torsion in its homology, this integral class is
well defined up to an integral coboundary so the same holds for
$c_\Fat^k$.
\end{proof}

The simplest example is $k=0$. Then
\[
    c_\Fat^0(\Gamma)=\sum_{v\in\Gamma^0}m(v)
    =-2\chi(\Gamma)=-2\chi(\Sig_\Gamma),
\]
ie, negative $2$ times the Euler characteristic of $\Gamma\simeq
\Sig_\Gamma$.

\subsection{Smooth families of fat graphs}

Suppose we have a smooth family of punctured surfaces, ie, a
smooth bundle $\Sig\to E\xrarrow{p}M$ where $M$ is a compact
smooth $n$--manifold with a fixed trivialization $E|\d
M=\Sig\times\d M$ over the boundary $\d M$ of $M$. If $\Sig$ is an
oriented surface of genus $g$ with $s\geq1$ boundary components
then $p\co E\to M$ is classified by a continuous mapping
\[
    f\co (M,\d M)\to (BM_g^s,\ast).
\]
By the simplicial approximation theorem we can choose any small
triangulation of $(M,\d M)$ and approximate $f$ by a simplicial
map
\[
    F\co (T(M),T(\d M))\to\N\bu\Fat_g^s
\]
where $\N\bu\Fat_g^s$ is the simplicial nerve of the category
$\Fat_g^s$ (fat graphs $\Gam$ whose surfaces $\Sig_\Gam$ have
genus $g$ wand $s$ boundary components). The following lemma
implies that $F$ can be chosen so that its image contains no fat
graphs of codimension $>n$.

\begin{lem}\label{FG:lem:link of codim c fat graphs}
Let $\Gam$ be any fat graph of codimension $c$. Then the full
subcategory of $\Fat/\Gam$ consisting of graphs of codimension
$<c$ is homotopy equivalent to a $(c-1)$--sphere.
\end{lem}

\begin{rem}\label{FG:rem:link of codim c fat graphs}
In any $k$ simplex $(\Gam_0\to\cdots\to\Gam_k)\in\N_k\Fat$, the last object $\Gam_k$ has the largest codimension. Therefore the subcategory of $\Fat/\Gam$ described in the lemma is the \emph{link} at $\Gam$ of the space of fat graphs of codimension $\geq c$.
\end{rem}

\begin{proof}
Suppose $v_1,\cdots,v_r$ are the vertices of $\Gam$ of $\codim>0$.
Let $c_i=\codim(v_i)$. Then $c=\sum c_i$. If $\Gam'$ is a fat
graph of codimension $c'<c$ which maps to $\Gam$ then in $\Gam'$
the vertices $v_1,\cdots,v_r$ must resolve into planar trees of
codimension $c_i'$ where $c_i'\lewq c_i$ and $\sum c_i'=c'<c$. In
other words, $\Gam'$ lies on the boundary of the product
\begin{equation}\label{FG:eq:product of associahedra}
  A_{c_1+3}(v_1)\times\cdots\times A_{c_r+3}(v_r)
\end{equation}
of the Stasheff polyhedra $A_{c_i+3}(v_i)$ associated with the
vertices $v_i$. But each Stasheff polyhedron is a disk so the
product (\ref{FG:eq:product of associahedra}) is a disk of
dimension $\sum c_i=c$ and $\Gam'$ lies on the boundary of this
disk, ie, it lies on a sphere of dimension $c-1$. (Actually,
$\Fat/\Gam$ is much larger since it contains infinitely many
isomorphic copies of each object so we get only a homotopy
equivalence with $S^{c-1}$.)
\end{proof}

\begin{prop}\label{FG:prop:no fat graphs of large codim}
We can choose the triangulation $(T(M),T(\d M))$ and the
simplicial map $F$ so that:
\begin{enumerate}
  \item The image of $F$ contains no fat graphs of codimension
  $>n$.
  \item If a vertex $v$ of $T(M)$ maps to a fat graph of
  codimension $n$ then the star of $v$ maps isomorphically to the
  product of Stasheff polyhedra (\ref{FG:eq:product of
  associahedra}).
\end{enumerate}
\end{prop}

\begin{proof}
By the lemma we may assume (1) and the condition that only
isolated vertices $v$ of $T(M)$ map to fat graphs of codimension
$n$. Then the link of such a vertex $v$ maps to the geometric
realization of the subcategory $\C$ of $\Fat/F(v)$ of graphs of
codimension $<n$ which is equivalent to an $n-1$ sphere by the
lemma. Consequently we have a well defined degree, say $d$. Now
modify the triangulation in the star of $v$ to include $d$ copies
of the $n$--disk (\ref{FG:eq:product of associahedra}). The
complement maps to $|\C|$.
\end{proof}

\subsection{The Witten cycle $W_{k}$}\label{subsec:Witten cycle}

Let $\W_{k}$ denote the full subcategory of $\Fat$ consisting of
fat graphs having a vertex of valence $\geq 2k+3$. Then it is
well-known that the realization
\[
    W_{k}=|\W_{k}|
\]
is a codimension $2k$ subset of $|\Fat|$ whose dual $[W_k]^\ast$
is a well-defined cohomology class for the mapping class group. In
the original finite dimensional model of
\cite{[Penner91:PoincareDual]} and \cite{[Kontsevich:Airy]},
$[W_k]^\ast$ is the Poincar\'{e} dual of the properly embedded
suborbifold $W_k$. For the category $\Fat$, we need to rely on
Proposition \ref{FG:prop:no fat graphs of large codim}. Given any
smooth family of punctured surfaces $\Sig\to E\to M$ as above
where $n=\dim M=2k$, we use Proposition \ref{FG:prop:no fat graphs
of large codim} to make $M$ ``transverse'' to $W_k$. A finite
number of vertices in the triangulation $T(M)$ will map to the
Witten cycle $W_k$. We count these with the orientation convention
of Definition \ref{FG:def:simplicial orientation of A2k+3} below.
To see directly that this gives a well-defined integer we need to
show that $[W_k]^\ast$ is zero on the link of a codimension $2k+1$
graph. Penner shows this geometrically in
\cite{[Penner91:PoincareDual]}. In this paper this will be a
consequence of our calculations since we prove that $[c_\Fat^k]$
is dual to the cycle $[W_{k}]$. The following theorem implies that
the signed intersection number of a $2k$ parameter family of fat
graphs with $W_k$ is given (up to a constant multiple) by
evaluation of the cocycle $c_\Fat^k$ on the family.

\begin{thm}\label{FG:thm:bk is dual to Wk}
The cohomology class $[c_\Fat^k]$ is a multiple of the dual
$[W_{k}]^\ast$ of the Witten cycle $W_{k}$. More precisely:
\begin{enumerate}
  \item The fat graph cocycle $c_\Fat^k$ vanishes on any
  $2k$--simplex of $\Fat$ which is disjoint from $\W_{k}$.
  \item If $\Gam$, $\Gam'$ are two objects of $\W_{k}$ of codimension $2k$
  (ie, having one vertex of valence $2k+3$ and all other
  vertices trivalent). Then $\Fat/\Gam\cong\Fat/\Gam'$ and $c_\Fat^k$ takes equal values on corresponding $2k$ simplices.
\end{enumerate}
\end{thm}

\begin{proof}
For (1) we note that the cyclic set cocycle $c_\Z^k$ is applied to
the angle sets of vertices. In order to be non-zero this set must
increase in size $2k$ times starting with $3$. The last vertex
must have valence $\geq 2k+3$ making the final fat graph an object of
$\W_{k}$.

For (2) we note that any object $\Gam''$ over $\Gam$ differs from $\Gam$ only at the vertex of valence $2k+3$. So, there is a corresponding object over $\Gam'$ equal to $\Gam''$ at this vertex and equal to $\Gam'$ away from the vertex. The cocycle $c_\Z^k$ is being applied only at this vertex. So, its value is the same on corresponding simplices.
\end{proof}

It remains to compute the value of the fat graph cocycle on the
category $\Fat/\Gam$. The realization of this category is the
Stasheff associahedron (if we choose one object from each
isomorphism class).

\subsection{The Stasheff associahedron}

It was Stasheff \cite{[Stasheff63:Polyhedron]} who constructed the
\emph{Stasheff associahedron} $A_n$ (also called the
\emph{Stasheff polyhedron}) and showed that it was an $n-3$ disk.
The reason it is called an associahedron is because it is an operad which continuously parametrizes the possible
ways to multiply, say, $n-1$ loops in any loop space. This
interpretation will not play a role in this paper.

For $n\geq3$ let $\A_n$ be the poset of all isomorphism classes of
planar trees with $n$ fixed leaves. This is a finite poset whose
elements have various interpretations which we use
interchangeably. One easy method to define this poset is to say
that its elements are sets $\Gam$ of unordered pairs $\{a,b\}$ of
distinct integers modulo $n$ so that:
\begin{enumerate}
  \item $b-a\nequiv\pm1 \mod n$.
  \item Whenever $a<b<c<d<a+n$, the two pairs $\{a,c\}$ and $\{b,d\}$
are not both in the set $\Gam$.
\end{enumerate}
In that case $\Gam$ represents the planar tree with vertices
$1,2,\cdots,n$ in cyclic order so that the two arcs $(i,i+1)$ and
$(j,j+1)$ bound two regions which come together at an internal
edge if and only if $\{i,j\}\in \Gam$. We will say that this edge
\emph{separates} the \emph{regions} $i$ and $j$.

The set $\A_n$ is ordered by inclusion. It has a unique minimal
element (the empty set in the above interpretation), so it is
obviously contractible. Also, the height of any element is at most
$n-3$. Thus the geometric realization $|\A_n|$ is $n-3$
dimensional.

We often view $\A_n$ as a category having a unique morphism $x\to
y$ iff $x\geq y$.

\begin{defn}\label{FG:def:Stasheff associahedron An}
The \emph{Stasheff associahedron} $A_n$ is defined to be the
geometric realization $|\A_n|$ of the category $\A_n$.
\end{defn}

Now consider the case $n=2k+3$. Let
\[
    \Gam_0>\Gam_1>\cdots>\Gam_{2k}
\]
be a nondegenerate $2k$--simplex in $\A_{2k+3}$. Then $\Gam_0$ is a
trivalent graph (interior vertices are trivalent) with exactly
$2k$ internal edges. The internal edges can be numbered
$e_1,e_2,\cdots,e_{2k}$ so that the edge $e_i$ collapses in the
$i^{\rm th}$ step (from $\Gam_{i-1}$ to $\Gam_i$).

Now number the internal vertices. The first edge $e_1$ has two
vertices which we label $v_0$, $v_1$ at random. For each other
edge $e_i$ let $v_i$ be the endpoint of $e_i$ which is the
furthest away from $v_0$. For $i\geq2$ this is also the endpoint
of $e_i$ furthest away from $v_1$ so the numbering of
$v_2,\cdots,v_{2k}$ remains unchanged if we switch $v_0,v_1$.

Let $a_1,a_2,a_3$ be, in cyclic order, the three regions which
come together at the vertex $v_0$ so that $e_1$ separates the
regions $a_1,a_3$. For $i=1,2,\cdots,2k$ let $b_i$ be the region
which touches the vertex $v_i$ without touching the interior of
the edge $e_i$.

\begin{lem}\label{FG:lem:labels are distinct} $a_1,a_2,a_3,b_1,\cdots,b_{2k}$ are distinct regions
(distinct elements of the set $\{0,1,2,\cdots,2k+2\}$.
\end{lem}

\begin{proof} For each region $r$ take the vertex in its boundary
which is closest to the vertex $v_0$. This sends $b_i$ to $v_i$
and $a_1,a_2,a_3$ to $v_0$.
\end{proof}

\begin{defn}\label{FG:def:simplicial orientation of A2k+3} The \emph{orientation} of the $2k$--simplex
$\Gam_0>\Gam_1>\cdots>\Gam_{2k}$ is defined to be the sign of the
permutation $(a_1,a_2,a_3,b_1,\cdots,b_{2k})$ of the set
$\{0,1,\cdots,2k+2\}$.
\end{defn}

Note first that this orientation is well-defined. If we reverse
(the names of) the vertices $v_0,v_1$ then $a_1,a_2,a_3,b_1$
become $a_3,b_1,a_1,a_2$ (an even permutation).

By Stasheff $A_{2k+3}=|\A_{2k+3}|$ is a $2k$--disk. Consequently,
the $2k$--simplices can be oriented so that their homological
boundary is the boundary sphere. We just need to check that we
have one of the two consistent orientation conventions.

\begin{prop}\label{FG:prop:orientation of A2k+3 is consistent}
Definition~\ref{FG:def:simplicial orientation of A2k+3} gives a
consistent orientation of the $2k$ simplices of the Stasheff
associahedron $A_{2k+3}$.
\end{prop}

\begin{proof} To check the consistency of our sign convention we only need
to show that the sign changes under both of the following
involutions:
\begin{enumerate}
  \item interchanging the numbering of two edges $e_i$ and
  $e_j$
  \item transforming the edge $e_1$ so that it goes the other way
  (separates $a_2,b_1$ instead of $a_1,a_3$).
\end{enumerate}
It is easy to see that both of these transformations changes the
sign of the permutation.

The first involution merely switches the labels $b_i$, $b_j$ (when
$i,j>1$). This is a transposition and thus odd. In the case $i=1$,
we first consider that case when $e_j$ is (geometrically) adjacent
to $e_1$. In that case we may assume that $v_0$ is the common
vertex of $e_1,e_j$. Then when we switch $e_1,e_j$ the distance
function to $v_0$ does not change so the vertices $v_i$ and the
regions $b_i$ do not change for $i\noteq 1,j$. The regions
$b_1,b_j$ are interchanged and the regions $a_1,a_2,a_3$ are
cyclically permuted. Consequently, the sign changes. All other
permutations can be obtained by composing these two operations.

The second transformation changes only the regions
$a_1,a_2,a_3,b_1$. They become $b_1,a_1,a_2,a_3$ (or
$a_2,a_3,b_1,a_1$), an odd permutation.

Finally, we need to explain why it suffices to show that (1) and
(2) switch the sign of the simplex.

Our $2k$--simplex has $2k+1$ faces. The last face does not count
since it is on the boundary sphere. The $0^{\rm th}$ face is given by
\[
    \Gam_1>\Gam_2>\cdots>\Gam_{2k}.
\]
In $\Gam_1$ the first edge $e_1$ is collapsed. This $2k-1$ simplex
is also the $0$--face of the $2k$--simplex
\[
    \Gam_0'>\Gam_1>\cdots>\Gam_{2k}
\]
where $\Gam_0'$ is obtained from $\Gam_0$ by operation (2) above.
Consequently, we need this operation to switch the sign of the
simplex so that the $0$--faces will homologically cancel.

The $i^{\rm th}$ face of the $2k$ simplex $\Gam_\ast$ is given by
deleting $\Gam_i$ from the sequence. However, when $0<i<2k$, this
is the same as deleting the $i^{\rm th}$ face of the $2k$--simplex
$\Gam_\ast'$ obtained from $\Gam_\ast$ by switching the labels
$e_i$ and $e_{i+1}$ (and their order of collapse). This is
operations (1).

Consequently, the boundary of the $2k$--chain given by the sum of
all $2k$--simplices with signs as above is equal to the sum of
their $2k$--faces forming the boundary $2k-1$ sphere.
\end{proof}

\subsection{Computation of $c_\Fat^k(W_{k})$}

We are now ready to compute the value of $c_\Fat^k$ on the
$2k$--chain given by the sum of all simplices in $\A_{2k+3}$ with
signs according to our convention. First we note that $c_\Fat^k=0$
on many of these $2k$ simplices. To have a chance to be non-zero we
need the graphs $\Gam_i$ in $\Gam_\ast$ to have the property that
they have one vertex of valence $i+3$ and all other internal
vertices to be trivalent. Let $B_k$ be the $2k$--chain consisting
of only these $2k$--simplices times appropriate signs:
\begin{equation}\label{FG:eq:Bk}
  B_k=\sum\sgn(a_1,a_2,a_3,b_1,\cdots,b_{2k})\Gam_\ast
\end{equation}
There are exactly
\begin{equation}\label{FG:eq:number of permutations}
  \frac{(2k+3)!}{6}
\end{equation}
permutations of the letters $a_1,a_2,a_3,b_1,\cdots,b_{2k}$ which
keep $a_1,a_2,a_3,b_1$ in cyclic order. These permutations
corresponds to the $2k$--simplex in the chain $B_k$ in a $2-1$
manner. In other words, $B_k$ has $(2k+3)!/3$ terms.

To see this it helps to go backwards. Start at the terminal graph
$\Gam_{2k}$. This has only one internal vertex. All regions
converge at that vertex. If $b_{2k}=j$ then the graph
$\Gam_{2k-1}$ is obtained from $\Gam_{2k}$ by pulling region $j$
away from the center of the graph leaving a vertex at which all
regions except for region $j=b_{2k}$. To get $\Gam_i$ from
$\Gam_{i+1}$ we pull away region $b_{i+1}$ from the central
vertex. This applies even in the case $i=0$ but in that case the
same result could have been obtained by ``pulling away'' the
region labelled $a_2$. In other words, if we switch $a_2,b_1$ and
we also switch $a_1,a_3$ then the resulting $2k$ simplex is the
same. There is a $2-1$ correspondence between permutations of the
letters $a_1,a_2,a_3,b_1,\cdots,b_{2k}$ which keep
$a_1,a_2,a_3,b_1$ in cyclic order and $2k$ simplices which occur
in the $2k$ chain $B_k$.

There is another way to look at this $2-1$ correspondence. The
formula for the combinatorial Miller--Morita--Mumford class applies
the cyclic set cocycle $c_\Z^k$ to exactly two vertices of the
graph $\Gam_0$ for each $2k$--simplex $\Gam_\ast$ which appears in
$B_k$. These are the vertices $v_0$, $v_1$ of $\Gam_0$ which merge
in $\Gam_1$. The permutations of the letters
$a_1,a_2,a_3,b_1,\cdots,b_{2k}$ which keep $a_1,a_2,a_3,b_1$ in
cyclic order are in natural $1-1$ correspondence with the
vertex-graph pairs upon which the cyclic set cocycles are to be
evaluated.

Now we consider the situation. We have a vertex $v_0$ in a graph
$\Gam_0$. The angles at $v_0$ are labelled $a_1,a_2,a_3$. As we
pass down to $\Gam_1,\Gam_2$, etc. the angles $b_1,b_2$, etc. are
added. Consequently, the value of the cyclic set cocycle $c_\Z^k$
is given by:
\begin{equation}\label{FG:eq:value of cZk on Gam v0}
  c_\Z^k(\Gam_\ast,v_0)=\frac{
    (-1)^kk!\sum_{i=1}^3\sgn(a_i,b_1,b_2,\cdots,b_{2k})
    }{
    (2k)!(2k+3)!/2
    }
\end{equation}
The value of the combinatorial cocycle $c_\Fat^k$ on the
$2k$--chain $B_k$ is therefore given by:
\begin{equation}\label{FG:eq:value of cFGk on Bk}
    c_\Fat^k(B_k)=
        \frac{(2k+3)!}{6}
\frac{
    (-1)^kk!E(X)
    }{
    (2k)!(2k+3)!/2
    }
    =\frac{
    (-1)^kk!E(X)
    }{
    3(2k)!
    }
\end{equation}
where $E(X)$ is the \emph{expected value} of the random variable:
\[
    X=\sgn(a_1,a_2,a_3,b_1,\cdots,b_{2k})\sum_{i=1}^3\sgn(a_i,b_1,b_2,\cdots,b_{2k})
\]
To compute $E(X)$ we note first that $X$ depends only on the
relative positions of $a_1,a_2,a_3$. For example, if they are
consecutive then $X=3$. More generally, if there are $p$ numbers
between $a_1$ and $a_2$, ie,
\[
    p=
  \begin{cases}
    a_2-a_1-1 & \text{if $a_2>a_1$}, \\
    a_2-a_1+2k+2 & \text{otherwise}.
  \end{cases}
\]
and $q$ numbers between $a_2$ and $a_3$, then
\begin{align*}
      \sgn(a_1,a_2,a_3,b_1,\cdots,b_{2k})
      \sgn(a_1,b_1,b_2,\cdots,b_{2k})& =(-1)^q \\
      \sgn(a_1,a_2,a_3,b_1,\cdots,b_{2k})
      \sgn(a_2,b_1,b_2,\cdots,b_{2k})& =(-1)^{p+q} \\
      \sgn(a_1,a_2,a_3,b_1,\cdots,b_{2k})
      \sgn(a_3,b_1,b_2,\cdots,b_{2k})& =(-1)^p
\end{align*}
so
\[
  X=(-1)^p+(-1)^q+(-1)^{p+q}=[1+(-1)^p][1+(-1)^q]-1.
\]
Let $Y$ be the random variable representing the number of spaces
between $b_1$ and $a_2$:
\[
    Y=
  \begin{cases}
    a_2-b_1-1 & \text{if $a_2>b_1$}, \\
    a_2-b_1+2k+2 & \text{otherwise}.
  \end{cases}
\]
Then $Y$ takes values in the set $\{1,2,\cdots,2k\}$ with
probabilities:
\begin{equation}\label{FG:eq:prob distribution of Y}
  P(Y=y)=\frac{y(2k+1-y)}{\sum_{j=1}^{2k}j(2k+1-j)}
\end{equation}
We leave it as an exercise for the reader to show that value of
the denominator is:
\begin{equation}\label{FG:eq:value of denominator in P(Y=y)}
  \sum_{j=1}^{2k}j(2k+1-j)=\frac13k(2k+1)(2k+2)
\end{equation}
(It suffices to verify this for $k=0,1,2,3$.)

Now we use the formula:
$
    E(X)=\sum_{y=1}^{2k}E(X|Y=y)P(Y=y).
$
We compute the conditional expected value of $X+1$:
\[
    E(X+1|Y=y)=\frac1{y(2k-y+1)}\sum_{p=0}^{y-1}\sum_{q=0}^{2k-y}
    [1+(-1)^p][1+(-1)^q]
\]
\[
    =\left(
\frac1{y}\sum_{p=0}^{y-1}[1+(-1)^p]
    \right)\left(
\frac1{2k-y+1}\sum_{q=0}^{2k-y}[1+(-1)^q]
    \right)
\]
If $y$ is odd then the right hand factor is $1$ and the left hand
factor is $1+\frac1{y}$. Similarly, if $y$ is even the left hand
fact or is $1$ and the right factor is $1+1/(2k-y+1)$.
Consequently,
\[
    E(X+1|Y=y)=
  \begin{cases}
    1+\frac1{y} & \text{if $y$ is odd}, \\
    1+\frac1{2k-y+1} & \text{if $y$ is even}.
  \end{cases}
\]
Using the apparent symmetry between $y$ and $2k-y+1$ we get:
\[
    E(X)=2\sum_{j=1}^kE(X|Y=2j)P(Y=2j)=\frac{2
\sum_{j=1}^k2j
    }{
\sum_{j=1}^{2k}j(2k+1-j)
    }
=\frac3{2k+1}
\]
Plugging this into (\ref{FG:eq:value of cFGk on Bk}) we get:
\begin{equation}\label{FG:eq:final value of cFGk on Bk}
        c_\Fat^k(B_k)=
\frac{
    (-1)^kk!E(X)
    }{
    3(2k)!
    }=
(-1)^k\frac{k!}{(2k+1)!}
\end{equation}
This proves:
\begin{thm}\label{FG:thm:bk is proportional to Wk}
\[
    [c_\Fat^k]=(-1)^k\frac{k!}{(2k+1)!}[W_{k}]^\ast
\]
as elements of $H^{2k}(\Fat;\QQ)\cong H^{2k}(M_g^s;\QQ)$.
\end{thm}

\subsection{Kontsevich cycles}

Suppose that $k_1,k_2,\cdots,k_r$ are positive integers. Then
Kontsevich \cite{[Kontsevich:Airy]} defined the cycles
\[
    W_{k_1,k_2,\cdots,k_r}
\]
in the space of fat graphs to be the set of all fat graphs having
$r$ vertices $v_1,\cdots,v_r$ with valences
$2k_1+3,2k_2+3,\cdots,2k_r+3$ respectively and no other vertices
of valence $>3$. Kontsevich also conjectured in
\cite{[Kontsevich:Airy]} that these were all related to the
Miller--Morita--Mumford classes. Arbanello and Cornalba
\cite{[Arbarello-Cornalba:96]} made this more precise by
conjecturing that their duals
\[
    [W_{k_1,k_2,\cdots,k_r}]^\ast\in
    H^{\deg(k_\ast)}(\Fat;\ZZ)=H^{\deg(k_\ast)}(M_g^s;\ZZ)
\]
should be polynomials in the adjusted Miller--Morita--Mumford
classes $\wt{\k}_j$. Here the \emph{degree} $\deg(k_\ast)$ of the
\emph{weight} $k_\ast=\{k_1,\cdots,k_r\}$ is given by
\[
    \deg(k_1,\cdots,k_r)=2k_1+\cdots+2k_r.
\]
The weights are partially ordered by $
\{j_1,\cdots,j_s\}<\{k_1,\cdots,k_r\}
$
if $s<r$ and there is an epimorphism of sets
$p\co \{1,\cdots,r\}\to\{1,\cdots,s\}$ so that
\[
    j_t=\sum_{i\in p^{-1}(t)}k_i.
\]
Because of the easy to understand nature of the combinatorial
classes $[c_\Fat^k]$ we get the following corollary of the
calculation in Theorem~\ref{FG:thm:bk is proportional to Wk}.

\begin{cor}\label{FG:cor:bk and Kontsevich cycles} Consider the
case when $k_1<k_2<\cdots<k_r$ occur with multipli\-cities
$n_1,\cdots,n_r$. Then
\[
    [c_\Fat^{k_1}]^{n_1}\cdots[c_\Fat^{k_r}]^{n_r}=
    \prod_{i=1}^{r}n_i!
    \left(
\frac{(-1)^{k_i}k_i!}{(2k_i+1)!}
    \right)^{n_i}
    \left[
W_{k_1^{n_1}\cdots k_r^{n_r}}
    \right]^\ast
    + \text{lower terms}
\]
where the ``lower terms'' refers to rational linear combinations
of dual Kontsevich classes $[W_{j_1\cdots j_m}]^\ast$ of weights
$\{j_1,\cdots,j_m\}$ less than
$\{(k_1)^{n_1},\cdots,(k_r)^{n_r}\}$ in the partial ordering
defined above.
\end{cor}

\begin{proof}
This is very elementary. Take the cup product of the fat graph
cocycles
\[
    c=c_\Fat^{p_1}\cup \cdots c_\Fat^{p_s}
\]
where $s=\sum n_i$ and $p_1,\cdots,p_s$ are $k_1,\cdots,k_r$ with
multiplicities $n_1,\cdots,n_r$.

When we evaluate the cocycle $c$ above on the $2d$ simplex
\[
    \Gam_\ast=(\Gam_0\to\cdots\Gam_{2d})
\]
where $2d=\deg(p_\ast)$ we will get zero unless there is a vertex
of $\Gam_0$ which increases in valence at each step until we reach
$\Gam_{2p_1}$, then we need a vertex of $\Gam_{2p_1}$ which
increases in valence until we reach $\Gam_{2p_1+2p_2}$, etc.
Assuming that $\Gam_{2d}$ has codimension $\leq2d$, this is
possible if and only if the graph $\Gam_{2d}$ lies in a Kontsevich
cycle $W_{j_1,\cdots,j_m}$ where
\[
    \{j_1,\cdots,j_m\}\leq\{p_1,\cdots,p_s\}
    =
    \{(k_1)^{n_1},\cdots,(k_r)^{n_r}\}.
\]

If $m=s$ and $j_i=p_i$ then the cocycle $c$ above, when evaluated
on $\Gam_\ast$ will be the product of coefficients
\[
    \frac{(-1)^{k_i}k_i!}{(2k_i+1)!}
\]
coming from each vertex of $\Gam_{2d}$ of valence $>3$ with a
factor of $\prod n_i!$ counting the number of permutations of
these vertices which preserve their valences.
\end{proof}

Corollary~\ref{FG:cor:bk and Kontsevich cycles} expresses
monomials in the $[c_\Fat^k]$'s as rational linear combinations of
the Kontsevich cycles in such a way that the change of basis
matrix is upper triangular with non-zero diagonal entries. The
inverse of this matrix is upper triangular with diagonal entries
inverse to those in Corollary~\ref{FG:cor:bk and Kontsevich
cycles}. Thus we get the following.

\begin{cor}\label{FG:cor:Kontsevich cycles in terms of bks}
\[
    \left[
W_{k_1^{n_1}\cdots k_r^{n_r}}
    \right]^\ast
    =
    \prod_{i=1}^r\frac1{
n_i!}\left( \frac{(2k_i+1)!}{(-1)^{k_i}k_i!}
    \right)^{n_i}
    [c_\Fat^{k_1}]^{n_1}\cdots[c_\Fat^{k_r}]^{n_r}
    + \text{lower terms}
\]
where the ``lower terms'' are rational linear combinations of cup
products of the form $[c_\Fat^{j_1}] \cdots [c_\Fat^{j_m}]$ with weights
$\{j_1,\cdots,j_m\}< \{(k_1)^{n_1},\cdots,(k_r)^{n_r}\}$. 
In particular,
these lower terms have smaller algebraic degree, ie, $m<\sum
n_i$.
\end{cor}

If we use the formula $[c_\Fat^k]=-2\k_k$ proved in Theorem \ref{Calc:thm:proportionality constant between MMM and b} below we obtain Corollary \ref{intro:cor:Kontsevich cycles in terms of MMM} in the introduction.

%
%

\section{Framed fat graphs}

In this section we will modify the arguments of the previous
section to prove the Witten conjecture, namely, the
\emph{adjusted} Miller--Morita--Mumford classes $\wt{\k}_k$ are dual
to $[W_{k}]$:
\[
    \wt{\k}_k=a_k[W_{k}]^\ast.
\]
The value of the coefficient $a_k$ is computed in the next
section. In other words, we will show in this section that the
analogue of Theorem~\ref{FG:thm:bk is dual to Wk} holds for
$\wt{\k}_k$.

The idea is to replace fat graphs with \emph{framed fat graphs}.
By the \emph{Framed Graph Theorem} \cite{[I:BookOne]} every graph
has a framed structure unique up to contractible choice. Thus
there is no obstruction to doing this. Framed fat graphs
$(\Gam,\f)$, although more complicated than the underlying fat
graphs $\Gam$, have the advantage in that they determine
\emph{framed functions} on the surface $\Sig_\Gam$. Consequently,
the adjusted Miller--Morita--Mumford class $\wt{\k}_k$ is given by
evaluating the (extended) cyclic set cocycle $c_{\Z^+}^k$ on the
$0$--cells of the framed structure.

By using a $D_{4k+4}$--equivariant framing for all trees in the
Stasheff associahedron $A_{2k+2}$ we conclude that the extended
cyclic set cocycle $c_{\Z^+}^k$ can only be non-zero in a
neighborhood of the Witten cycle $W_{k}$, just as in the case of
fat graphs.

\subsection{The Framed Function Theorem}

Before we explain the Framed Graph Theorem we need to review the
definition and statement of the Framed Function Theorem
\cite{[I:FF]}. First we recall that a \emph{generalized Morse
function} (GMF) on a compact smooth manifold $M$ is a smooth
function $f\co M\to\RR
$
which has only Morse and birth-death
singularities. \emph{Birth-death points} are the unique
codimension $1$ singularities of smooth functions. They are points
at which the function $f$ can be written as
\[
    f(x)=x_0^3+\sum_{i\noteq0}\pm x_i^2+C
\]
with respect to some local coordinate system. Here $C$ is a
constant and the number of negative signs in the sum $\sum\pm
x_i^2$ is called the \emph{index} of the birth-death point.

A \emph{generic unfolding} of a birth-death point is given by
\[
    f_t(x)=x_0^3+tx_0+\sum_{i\noteq0}\pm x_i^2+C
\]
for $t$ in a small interval $(-\e,\e)$ about $0$. Note that when
$t>0$ the function $f_t$ has no critical points (in the coordinate
neighborhood). For $t<0$ the function $f_t$ has two critical
points in ``cancelling position.''

\begin{defn}\label{MMMW:def:framed functions}
A \emph{framed function} on a smooth manifold $M$ is a GMF
$f\co M{\to}$ $\RR$ together with a tangential framing of the nonpositive
eigenspace of the second derivative $D^2f$ at each critical point
so that the last framing vector points in the positive cubic
direction ($\frac{\d}{\d x_0}$) at every birth-death point. In a
family of framed functions we assume that the framing vectors vary
continuously.
\end{defn}

\begin{thm}\label{MMMW:thm:framed function theorem}{\rm\cite{[I:FF]}}\qua
The space of framed functions on $M$ is $\dim M-1$ connected.
\end{thm}

This theorem is not quite good enough in the present situation. We
need to know that the space of framed functions on a compact
surface is contractible. We cannot prove this but we have another
theorem which is just as good for our purposes: the Framed Graph
Theorem which implies that compact surfaces admit canonical framed
functions which are unique up to contractible choice. These framed
functions have the property that they have critical points only in
indices $0,1$.

\subsection{The Framed Graph Theorem}

We review the concept of a framed graph. It helps to keep in mind
that this is the combinatorial structure associated with a framed
function. These framed functions will have only three kinds of
singularities:
\begin{enumerate}
  \item Morse points of index $0$ (local minima).
  \item Morse points of index $1$ (maxima in dim 1 and
  saddle points in dim 2).
  \item Birth-death points of index 0 (where the above two
  cancel).
\end{enumerate}
We note that the framing of the framed function gives an
\emph{orientation} of the index $1$ Morse points resulting in an
oriented $1$--dimensional cell complex as the core of our manifold
(the domain of the framed function).

\begin{defn}\label{MMMW:def:framed graphs}
A \emph{framed graph} is an oriented $1$--dimensional cell complex
$X$ in which $1$--cells are allowed to be attached on the interiors
of other $1$--cells in some partial ordering of the set of cells
and which have a designated collection of \emph{collapsing pairs}.
Collapsing pairs consist of a $1$--cell and target $0$--cell.
\end{defn}

\begin{figure}[ht!]
\cl{\includegraphics[width=1.8in]{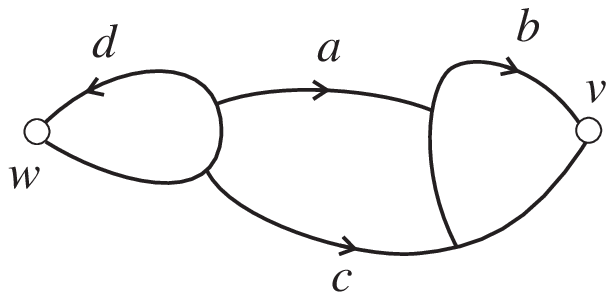}}
\nocolon\caption{}\label{MMMW:fig:02}
\end{figure}

To expand further on our analogy with framed functions we note
that in a $1$--parameter family of framed functions, the cells
corresponding to the critical points will slide over each other
discontinuously. In the combinatorial setting we need to do this
continuously. Thus the $1$--cells need to be allowed to attach to
the interiors of other $1$--cells. For functions we note that the
critical points are ordered according to critical value and that
cells can only slide over cells lower in this ordering. Also a
birth-death point can be interpreted as a pair of cancelling
critical points with the same critical value. Thus the partial
ordering of the cells in our framed graph should make collapsing
pairs equivalent in the partial ordering so that, eg, if $(e,v)$
forms a collapsing pair and $e'$ is another $1$--cell, if $e'$
attaches to $v$ then $e$ cannot be attached to the interior of
$e'$ (since $e'<e\sim v<e'$ would be a contradiction). Thus in
Figure~\ref{MMMW:fig:02} neither $(b,v)$ nor $(d,w)$ is allowed to
be a collapsing pair.

In order to understand the morphisms between framed graphs we need
to look at the \emph{underlying graph} $\Gam$ of a framed graph
$X$. This graph will have vertices of valence $2$ corresponding to
$0$--cells which occur in the centers of edges. The leaves
(vertices of valence $1$) are also necessarily $0$--cells. However,
vertices of valence $\geq3$ may or may not be $0$--cells.

If a vertex of the underlying graph $\Gam$ is not a $0$--cell then
it must lie on the interior of a $1$--cell. Consequently, two of
the incident half-edges must form part of a $1$--cell and the other
incident half-edges must be ends of other $1$--cells. (A $1$--cell
cannot be attached to the interior of itself by the partial
ordering rule and a $1$--cell in a cancelling pair cannot have
source equal to its target $0$--cell by the same rule.)

Strictly speaking, morphisms of framed graphs will be graph
epimorphisms as in the category of fat graphs. However, we want to
interpret them as \emph{handle slides}. For example, if a $1$--cell
$a$ attached to the center of another $1$--cell $b$ then by
collapsing one or the other half of $b$ we obtain what looks like
a sliding of the attaching map of $a$ to one or the other endpoint
of $b$. Thus a collapse of a portion of a $1$--cell can be viewed
as a handle slide.

Collapsing pairs are allowed to collapse to their source. They are
also allowed to ``unpair.'' This is analogous to the possible
deformations of a birth-death point of index $0$. It can either
resolve into a nonsingular function (collapsing the cell pair) or
it can resolve into a pair of distinct Morse points. The latter
possibility we represent as ``unpairing'' the cancelling pair,
ie, keeping the pair but removing their designation (demoting
them).

To emphasize the fact that we consider a framed graph $X$ as being
an embellishment of the underlying combinatorial graph $\Gam$ we
write $X=(\Gam,\f)$ where $\f$ represents the \emph{framed
structure} of $\Gam$ given by $X$.

\begin{defn}\label{MMMW:def:morphism of framed graphs}
A \emph{morphism} of framed graphs $(\Gam_0,\f)\to(\Gam_1,\psi)$
is a graph epimorphism $f\co \Gam_0\to\Gam_1$ with the property that:
\begin{enumerate}
  \item The collapsing trees $f^{-1}(v)$ do not contain any entire
  $1$--cells unless they form part of a collapsing pair in
  $\Gam_0$.
  \item Every $1$--cell of $\Gam_1$ is the image of a $1$--cell of
  $\Gam_0$ with the induced orientation.
  \item Every $0$--cell of $\Gam_1$ is the image of a $0$--cell of
  $\Gam_0$.
  \item Every collapsing pair of $\Gam_1$ is the image of a
  collapsing pair of $\Gam_0$.
\end{enumerate}
\end{defn}

We leave it to the reader to show that morphisms are closed under
composition. Also we note that condition (2) implies condition
(3). Another very useful observation is that $\Gam_0\simeq\Gam_1$
and therefore they have the same Euler characteristic. A certain
number of collapsing pairs will disappear according to (1), but,
in order to preserve the Euler characteristic, the remaining
$0$--cells of $\Gam_0$ must map to distinct $0$--cells of $\Gam_1$.
This implies the following.

\begin{prop}\label{MMMW:prop:what happens to collapsing pairs}
Under a morphism $f\co (\Gam_0,\f)\to(\Gam_1,\psi)$ each collapsing
pair $(e,v)$ of $\Gam_0$ will be transformed in one of three ways:
\begin{enumerate}
  \item $(e,v)$ will collapse to the source of $e$. (This removes
  one $0$--cell and one $1$--cell from the framed graph.)
  \item $(e,v)$ maps to a collapsing pair $(f(e),f(v))$ in
  $\Gam_1$.
  \item $(e,v)$ is \emph{unpaired}, ie, $f(e)$ is a $1$--cell and
  $f(v)$ is a $0$--cell but $f(e),f(v)$ do not form a collapsing
  pair.
\end{enumerate}
\end{prop}

For the purpose of giving the definition of a morphism
(\ref{MMMW:def:morphism of framed graphs}) we allowed our
underlying graph to have vertices of valence $2$. Normally this is
not useful. We define the \emph{reduced} graph of a framed graph
to be the underlying graph with the valence $2$ vertices
``smoothed'' in the sense that the two incident half edges are
joined together. This causes a problem only if the graph is a
circle. We regard the placement of the $0$--cells of valence $2$ to
be part of the framed structure $\f$ of the reduced graph.

The \emph{Framed Graph Theorem} says that any family of graphs
admits a family of framed structures and that this structure can
be specified on any closed subcomplex of the parameter space.
(Consequently, the framed structure is unique up to framed
homotopy.) However, in this paper we have not defined the notion
of ``family of graphs'' for arbitrary graphs. We have however
talked about families of fat graphs.

\subsection{Framed fat graphs}

A \emph{framed fat graph} is simply a fat graph with a framed
structure. Given a morphism of framed graphs
$(\Gam_0,\f)\to(\Gam_1,\psi)$ and a fat graph structure on
$\Gam_0$ there is an induced fat graph structure on $\Gam_1$ since
every tree in a subdivided $\Gam_0$ induces a cyclic ordering of
incident edges. A \emph{morphism} of framed fat graphs is a
morphism of framed graphs with compatible fat graph structures.

The Framed Graph Theorem implies that any family of fat graphs
(ie, continuous mapping into the geometric realization $|\Fat|$)
admits a continuous family of framed structures which can be
specified on any subcomplex of the parameter space.

The main theorem is that the adjusted Miller--Morita--Mumford class
$\wt{\k}_k$ is given by applying the extended cyclic set cocycle
$c_{\Z^+}^k$ to the $0$--cells of the framed structure. This
follows from an examination of the framed functions on the
punctured surfaces $\Sig_\Gam$ which are associated to the framed
fat graphs.

\subsection{Relation to Morse theory}\label{subsec:Morse}

Suppose that $\Sig\to E\xrarrow{p} M$ is a smooth oriented surface bundle
with a fiberwise Riemannian metric. In case $M$ has a boundary $\d M$, we
assume that $E|\d M=p^{-1}(\d M)$ is a trivial bundle $E|\d M\cong \d
M\times\Sig$. Let $\d_v E$ be the subbundle
of $E$ with fiber $\d\Sig$ (so that $\d E=\d_vE\cup E|\d M$). Since $\d\Sig$ is a disjoint union of circles the mapping $\d_vE\to M$ is the composition of an oriented circle bundle $\d_vE\to \wt{M}$ with a covering map $\pi\co \wt{M}\to M$.
We recall the definition of the Miller--Morita--Mumford class. Take the Euler class $e_E\in H^2(E,\d_vE;\ZZ)$ of the vertical tangent bundle of $E$. Then the \emph{adjusted Miller--Morita--Mumford class} is given by $\wt\k_k=p_\ast e_E^{k+1}\in H^{2k}(M,\d M;\ZZ)$ where
\[
	p_\ast\co H^\ast(E,\d E;\ZZ)\to H^{\ast-2}(M,\d M;\ZZ)
\]
is the \emph{push-down} operator. See \cite{[MoritaSurvey]} for more details. The usual Miller--Morita--Mumford class $\k_k\in H^{2k}(M,\d M;\ZZ)$ is related to $\wt\k_k$ by
\[
	\k_k=\wt\k_k+\pi_\ast(\g^k)
\]
where $\g\in H^2(\wt{M},\d\wt{M};\ZZ)$ is the Euler class of the oriented circle bundle $\d_vE\to \wt{M}$ and $\pi_\ast\co H^\ast(\wt{M},\d\wt{M};\ZZ)\to H^\ast(M,\d M;\ZZ)$ is the push-down map.

Let $f\co E\to \RR$ be any generic smooth function whose fiberwise
gradient $\grad f_t$ points outward along the boundary $\d\Sig_t$
of every fiber $\Sig_t$.  Since $E|\d M\cong \Sig\times\d M$ we may assume that $f_t\co \Sig_t=\Sig\to\RR$ is a fixed Morse function for all $t\in\d M$.

Since $f$ is generic, its vertical singular set $\Sig(f)$ is a
codimension $2$ submanifold of $E$ which is disjoint from $\d_vE$
and which is a product bundle over $\d M$. Since $\Sig(f)$ is the
inverse image under the vertical gradient $\grad f_t$ of the zero
section of the vertical tangent bundle of $E$ we have the
following well-known observation.

\begin{prop}\label{WWWM:prop:the singular set is dual to the Euler
class} The vertical singular set $\Sig(f)$ is Poincar\'{e} dual to
the Euler class $e_E$ of the vertical tangent bundle of $E$.
\end{prop}

\begin{rem} This is not quite correctly stated since $\Sig(f)$
meets the boundary of $E$. What we mean is that cup product with
the Euler class $e_E$ is equivalent to restriction to $\Sig(f)$.
In particular, the $k+1$st power $e_E^{k+1}$ is equivalent to the
restriction to $\Sig(f)$ of $e_E^k$ in the sense that they have
the same push-down in $H^{2k}(M,\d M;\RR)$.\end{rem}

Now suppose that $f\co E\to \RR$ is a fiberwise framed function which
is still fixed over $\d M$. The framed structure gives a
trivialization of the vertical tangent bundle along birth-death
points and saddle points (and maxima if there are any). Also we
can choose a trivialization over $\d M$. This gives the following.

\begin{cor}\label{WWWM:cor:MMM classes are given by 0 cells}
The adjusted Miller--Morita--Mumford classes $\wt{\k}_k$ are given by pushing down the restriction to
the index $0$ singular set $\Sig^0(f)$ of the $k^{\rm th}$ power $e_E^k$
of the vertical Euler class $e_E$:
\[
    \wt{\k}_k=p_\ast(e_E^k|\Sig^0(f))
\]
\end{cor}

On the category of framed fat graphs we need to use the extended
cyclic set cocycle $c_{\Z^+}^k$ on the set of $0$--cells. This
takes into account the trivialization of the vertical tangent
bundle at the places where the $0$--cell is part of a collapsing
pair.

\begin{defn}\label{MMMW:def:alpha k on frFat}
Suppose that $\Gam_\ast=(\Gam_0\to\Gam_1\to\cdots\to\Gam_{2k})
$
is a $2k$--simplex in the category of framed fat graphs. Then let
$z^k(\Gam_\ast)$ be given by
\[
    z^k(\Gam_\ast)=\sum c_{\Z^+}^k(C(v)\to C(f_1(v))\to\cdots\to
    C(f_{2k}(v)))
\]
where the sum is over all $0$--cells $v$ of $\Gam_0$ which survive
(do not collapse) in $\Gam_{2k}$. We call this the \emph{$0$--cell
cocycle}.
\end{defn}

Since a framed structure on a fat graph $\Gam$ gives a framed
function on the punctured surface $\Sig_\Gam$ (See the last
section of \cite{[I:BookOne]}), Corollary~\ref{WWWM:cor:MMM
classes are given by 0 cells} gives the following.

\begin{cor}\label{MMMW:cor:MMM is given by 0 cell cocycle}
The adjusted Miller--Morita--Mumford classes $\wt{\k}_k$ are given
by the $0$--cell cocycle applied to the transverse families of
framed fat graphs associated to the surface bundle $\Sig\to E\to
M$.
\end{cor}

By \emph{transverse} we mean that the underlying family of fat
graphs satisfies Proposition~\ref{FG:prop:no fat graphs of large
codim}.


\subsection{Localizing to Witten cycles}

We would like to show that the only terms in $z^k(\Gam_\ast)$
which can be non-zero are the ones in which the last vertex
$f_{2k}(v)$ has valence $2k+3$, ie, $\Gam_{2k}$ is in the Witten
cycle $W_{k}$. This is very close to being true.

There are two ways that this might fail. The first is that framed
fat graphs have $0$--cells of valence $2$. The second is that the
extended cyclic set cocycle $c_{\Z^=}^k$ can be non-zero even if
the last two sets have the same size. Thus the worst case is:
\begin{equation}\label{WWWM:eq:worse case is 2,3,...,2k-1,2k-1}
  |C(v)|=2\ ,\ |C(f_1(v))|=3\ ,\ \cdots\ ,\ |C(f_{2k-1}(v))|=2k+1=|C(f_{2k}(v))|
\end{equation}
where $v,f_1(v),\cdots,f_{2k-1}(v)$ are paired $0$--cells and
$f_{2k}(v)$ is unpaired.

In order to prevent this we will show that the framed structure on
a transverse $2k$--parameter family of fat graphs can be chosen so
that all ``massive'' points are unpaired. By a \emph{massive
point} we mean a $0$--cell of valence $\geq k+3$ (codim $\geq k$).

If $k\geq2$ then $2k+1\geq k+3$ so the worse case
(\ref{WWWM:eq:worse case is 2,3,...,2k-1,2k-1}) will not occur if
all massive points are unpaired. The case $k=1$ is treated
separately.

Assume for a moment that $k\geq2$ and all massive points are
unpaired. Then the vertex $f_{2k-1}(v)$ must be unpaired and the
worst case (\ref{WWWM:eq:worse case is 2,3,...,2k-1,2k-1}) does
not occur. Instead the worst case is
\begin{equation}\label{MMMW:eq:worst case no 2}
  |C(v)|=2\ ,\ |C(f_1(v))|=3\ ,\ \cdots\ ,\ |C(f_{2k}(v))|=2k+2
\end{equation}
with all $0$--cells $f_i(v)$ unpaired. However, this is not
possible since the associahedron $A_{2k+2}$ giving the unfolding
of a vertex of valence $2k+2$ is a $2k-1$ disk which forces the
$2k$--simplex (\ref{MMMW:eq:worst case no 2}) to be degenerate. To
make this argument valid we need to choose a
$D_{4k+4}$--equivariant framing of the trees in the Stasheff
associahedron $A_{2k+2}$. We will give an explicit such framing
which is compatible with our first condition, (ie, massive
points are unpaired).

Of course, any triangulation of a $2k-1$ disk will have no
nondegenerate $2k$--simplices. The essential point is the dihedral
($D_{4k+4}$) symmetry of the framing which allows us to give a
product structure to the framed structures in a neighborhood of
the arcs of $W_{k-\frac12}$ graphs (having a $2k+2$ valent vertex)
in our generic $2k$--parameter family of fat graphs. The final
result will be that the $0$--cell cocycle $z^k$ will only be
non-zero at two places:
\begin{enumerate}
  \item The final vertex $f_{2k}(v)$ has valence $2k+3$ forcing the final
fat graph $\Gam_{2k}$ to lie in the Witten cycle $W_{k}$.
  \item The final vertex has valence $2k+2$ and lies near the
  boundary of the associahedron $A_{2k+3}$ where the dihedral
  symmetry of the lower associahedron $A_{2k+2}$ is broken.
\end{enumerate}
This is explained in detail in the last section.

\subsection{Massive points}

The first step is to arrange for all massive vertices to be
unpaired. To do this we need the following relative version of the
Framed Graph Theorem. Fortunately, it follows from the absolute
version.

Let $\Gam\co \nP\to\Fat$ be a functor from any small category $\nP$
into the category of connected fat graphs. Then we get an induced
continuous family of fat graphs $\Gam(t),t\in|\nP|$. Suppose that
each $\Gam(X),X\in\nP$ has two disjoint subgraphs $A(X)$ and
$B(X)$ which contain all the vertices of $\Gam(X)$. Thus
\[
    E(X)=\Gam(X)-(A(X)\cup B(X))
\]
is a disjoint union of open edges. Suppose that this decomposition
is natural in the sense that every morphism $X\to Y$ in $\nP$
sends $A(X)$ to $A(Y)$, $B(X)$ to $B(Y)$ and $E(X)$ to $E(Y)$.
Then we get continuous families of subgraphs
$A(t),B(t)\subseteq\Gam(t)$ for $t\in|\nP|$ and the number of
elements of $E(t)$ is locally constant on $|\nP|$.

More generally, suppose that $K$ is a closed subcomplex of (a
subdivision of) $|\nP|$. For each $t\in K$ suppose we have
subgraphs $A(t),B(t)$ of $\Gam(t)$ containing all the vertices of
$\Gam(t)$ so that
\[
    E(t)=\Gam(t)-(A(t)\cup B(t))
\]
is a disjoint union of open edges which vary continuously with
$t\in K$ and the number of which is locally constant over $K$.

\begin{thm}\label{WWWM:thm:extended framed graph theorem}
Let $A(t),B(t)\subseteq \Gam(t)$, $t\in K$, be as described above.
Then any continuous family of framed structures on the subgraphs
$A(t)$, $t\in K$, extends to a continuous family of framed
structures on $\Gam(t)$ with the property that:
\begin{enumerate}
  \item $B(t)$ is a subgraph of $\Gam(t)$ for all $t\in K$.
  \item None of the $1$--cells in $E(t)$ forms a collapsing pair
  with a vertex in $A(t)\cup B(t)$.
\end{enumerate}
Furthermore, the relative version of this statement is true: The
framed structure on $\Gam(t)$ extending the given framed structure
on $A(t)$ can be specified for $t\in L$ for any subcomplex $L$ of
$K$ provided that the above two conditions hold over $L$.
\end{thm}

\begin{proof}
The absolute version is easy. Using the Framed Graph Theorem we
can simply choose a continuous framed structure on $B(t)$ for all
$t\in K$ and then take the following framed structure on the edges
of $E(t)$:
\[
    \longrightarrow 0\longleftarrow
\]
(Place a $0$--cell in the center and form two $1$--cells oriented
towards the center.)

For the relative case we use the relative version of the Framed
Graph Theorem to construct a continuous family of framed
structures on $B(t)$ extending the given structure over $L$. Then
we use the fact that the space of framed structures on an interval
with specified boundary behavior is contractible (the case
$\FG_{11}(I)$ of Theorem 8.3.2 of \cite{[I:BookOne]}).
\end{proof}

Now let $\nP_{2k}$ denote the full subcategory of $\Fat$
consisting of fat graphs of codimension $\leq 2k$ with at least
one massive vertex, ie, with valence $\geq k+3$. We want to
construct subgraphs $A(t)$ of $\Gam(t)$ for each $t\in|\nP_{2k}|$
containing the massive vertices of $\Gam(t)$ and we want to use
Theorem~\ref{WWWM:thm:extended framed graph theorem} above to
specify the framed structure on $A(t)$ (so that the massive vertex
is not a paired $0$--cell).

Consider a $p$--simplex $\Gam_0\to\cdots\to\Gam_p
$
in $\nP_{2k}$. If $\Gam_0$ has two massive vertices then it has
codimension $2k$ and $\Gam_i\cong\Gam_0$ for all $i$. So suppose
that $\Gam_0$ has a unique massive vertex $v$. Let $f_i(v)\in
\Gam_i$ be the image of $v$.

Let $T_i\subseteq\Gam_i$ be the inverse image of
$f_p(v)\in\Gam_p$. Then $T_i$ determines a continuous family of
subgraphs $T(t)\subseteq\Gam(t)$ for all $t\in \Delta^p$. However,
this is not a suitable choice for $A(t)$ since these are graphs
having vertices connected to a massive vertex by two edges.

Our second attempt is to define $A(t)$ to be the subgraph of
$\Gam(t)$ consisting of the massive vertex and any connected
``short edges'' where by a \emph{short edge} we mean an edge of
length $\leq\ll$ for some fixed $\ll\approx \frac12$ and by
\emph{connected} we mean either directly or through other short
edges. (Thus $A(t)$ is connected with all edges of length
$\leq\ll$.) However, there is one problem with this definition.
There might be edges of length slightly larger than $\ll$ whose
length decreases to $\LL$. We will avoid this and other unseen
problems with a more formal approach.

First, choose a small positive real number $\e$.
($0<\e<\frac1{2k+3}$ is sufficient.) Then, for each $k\leq
j\leq2k$ and $t\in|\nP_{2k}|$ we define:
\begin{enumerate}
  \item The \emph{closed $j$--weight} of $\Gam(t)$ to be the
  codimension of the massive point of $\Gam(t)$ plus the number of
  edges of length $\leq j\e+\frac12$ (call them \emph{$j$--short edges})
  which are connected to the massive point either directly or through other $j$--short edges.
  \item The \emph{open $j$--weight} of $\Gam(t)$ is defined to be the
  codimension of the massive point plus the number of edges of
  length strictly less than $j\e+\frac12$ (call them \emph{strictly $j$--short
  edges}) which are connected to the massive point by strictly
  $j$--short edges.
\end{enumerate}

Note that closed $j$--weight is upper semicontinous (being bounded
above is an open condition) and open $j$--weight is lower
semicontinuous (being bounded below is open). Also
\[
    \text{open $j$--weight $\leq$ closed $j$--weight}.
\]
For each $0\leq j\leq k$ let $K_j$ be the set of all
$t\in|\nP_{2k}|$ so that the closed $j$--weight of $\Gam(t)$ is
$\geq k+j$ and equal to the open $(j+1)$--weight of $\Gam(t)$. This
is a closed condition since it can be written as:
\[
    \text{open $(j+1)$--weight$(\Gam(t))\ -$ closed $j$--weight$(\Gam(t))\leq0$}
\]

\begin{lem}\label{MMMW:lem:|P|=cup Kj}
$|\nP_{2k}|$ is a union of the closed subcomplexes $K_j$ for
$0\leq j\leq k$.
\end{lem}

\begin{proof}
For any $t\in|\nP_{2k}|$ let $\ll_t\co I\to\NN $ be defined as
follows. $\ll_t(x)$ is equal to the codimension of the massive
point of $\Gam(t)$ plus the number of attached short edges where
\emph{short} means of length $\leq x$. Then $K_j$ is the set of
all $t\in|\nP_{2k}|$ so that $\ll_t(j\e+\frac12)\geq k+j$ and
$\ll_t$ is constant on the half-open interval $
[j\e+\tfrac12,(j+1)\e+\tfrac12)$.

For any $t\in|\nP_{2k}|$ we note that $\ll_t(\frac12)\geq k$.
Consequently, there exists a $j\geq 0$ so that
$
    \ll_t(j\e+\tfrac12)\geq k+j.
$ Choose the largest such $j\leq k$. If $t\notin K_j$ then we must
have
$
    \ll_t(x)>\ll_t(j\e+\tfrac12)
$
for some $x<(j+1)\e+\frac12$. This implies that
$
    \ll_t((j+1)\e+\tfrac12)\geq k+j+1
$
contradicting the maximality of $j$.
\end{proof}

For each $t\in K_j$, let $A_j(t)\subseteq\Gam(t)$ be the union of
the massive points and all connected $j$--short edges (of length
$\leq j\e+\frac12$). Let $E_j(t)$ be the set of edges connected to
$A_j(t)$. By definition of $K_j$, each edge in $E_j(t)$ has length
$\geq(j+1)\e+\frac12$. Consequently, the number of edges in
$E_j(t)$ is a locally constant function of $t\in K_j$ so the
conditions of Theorem~\ref{WWWM:thm:extended framed graph theorem}
are satisfied.

Now we can specify the framing of $A_j(t)$. We take the $0$--cells
to be the vertices and the $1$--cells to be the edges oriented away
from the massive point. We designate a $1$--cell and its target
vertex to be a collapsing pair if and only if its length is
$\leq\frac14$. We call this framing the \emph{radial framing} of
$A_j(t)$ centered at the massive point.

We extend this framing to the rest of $\Gam(t)$ using
Theorem~\ref{WWWM:thm:extended framed graph theorem} by downward
induction on $j$. To check the relative condition note that if
$t\in K_\l$ for $\l>j$ then $A_j(t)\subseteq A_\l(t)$ so $A_j(t)$
is already radially framed. (And any edges of $A_\l(t)$ which are
not in $A_j(t)$ have length $>\frac12$ so are not paired.) By the
framed graph theorem this framing can also be extended to all fat
graphs of codimension $\leq2k$. This shows:

\begin{lem}\label{MMMW:lem:massive points can be unpaired}
Framed structures can be chosen for all fat graphs of codimension
$\leq2k$ so that all massive points are unpaired.
\end{lem}

\subsection{The case $k=1$}

Now consider the special case $k=1$. The argument using massive
points fails in this case. However, there is a simple reason that
the desired statement still holds, ie, the $0$--cell cocycle
$z^1(\Gam_\ast)$ is still zero when $f_1(v)$ is a paired $0$--cell.
The reason is that, when $v$ is a paired $0$--cell of valence $2$
and $f_1(v)$ is a paired $0$--cell of valence $3$, there is only
one possible geometry as shown in Figure~\ref{MMMW:fig:03}. If
$f_2(v)$ is unpaired with valence $3$ then:
\[
    c_{\Z^+}^1(C(v)\to C(f_1(v))\to C(f_2(v)))=-s_1(C(v)\to
    C(f_1(v)))=\frac{\sgn(aba)}{2!\cdot2\cdot3}=0
\]
Therefore:

\begin{figure}[ht!]
\cl{\includegraphics[width=2.7in]{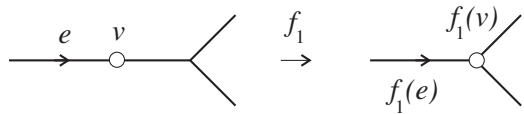}}
\caption{Unique geometry allowing $f_1(v)$ to be paired with
valence $3$}\label{MMMW:fig:03}
\end{figure}

\begin{lem}\label{MMMW:lem:the last vertex has codim 2k-1}
For all $k\geq1$ there are no non-zero terms in the $0$--cell
cocycle unless the last vertex $f_{2k}(v)$ has codimension $\geq
2k-1$ (assuming that massive points are unpaired).
\end{lem}

\subsection{Equivariant framing of the Stasheff associahedron $A_{2k+2}$}

An equivariant framing of the Stasheff associahedron $A_n$ for any
$n$ is easy to describe. We simply take the radial framing
centered at the center of mass of each object of $\A_n$
(Definition~\ref{FG:def:Stasheff associahedron An}).

\begin{figure}[ht!]
\cl{\includegraphics[width=1.5in]{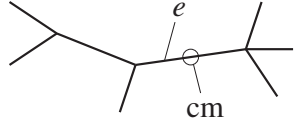}}
\caption{Center of mass (circle) for tree with 6
leaves}\label{MMMW:fig:500}
\end{figure}

\begin{defn}\label{MMMW:def:center of mass}
Suppose that $\Gam$ is a tree with $n$ leaves and no bivalent
vertices (equivalent to an object of $\A_n$). Then the
\emph{center of mass} of $\Gam$ is defined to be either:
\begin{enumerate}
  \item The midpoint of the unique edge $e$ of $\Gam$ with the
  property that each endpoint of $e$ is connected to half the
  leaves of $\Gam$ by paths disjoint from the interior of $e$ (Figure~\ref{MMMW:fig:500}).
  \item If no such edge exists then the unique vertex $v$ of
  $\Gam$ having the property that no component of $\Gam-e$ has
  more than half the leaves of $\Gam$.
\end{enumerate}
\end{defn}

We leave it to the reader to verify that the center of mass is
well-defined. Note that case (1) is possible only if $n$ is even.

\begin{prop}\label{MMMW:prop:properties of center of mass}
Let $n=2k+2$ or $2k+3$.
\begin{enumerate}
  \item Any morphism $\Gam_1\to\Gam_2$ in $\A_n$ sends the center
  of mass of $\Gam_1$ to the center of mass of $\Gam_2$.
  \item If $\Gam\in \A_n$ has a massive vertex $v$ then $v$ is the
  center of mass.
\end{enumerate}
\end{prop}

\begin{figure}[ht!]
\cl{\includegraphics[width=2.2in]{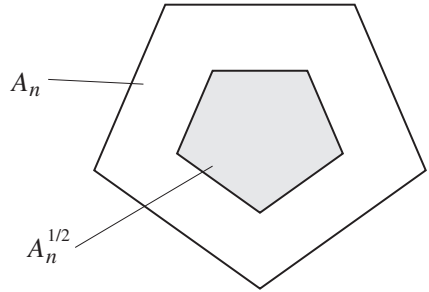}}
\caption{The small copy $A_n^{1/2}$ does not touch any other
$A_n$.}\label{MMMW:fig:501}
\end{figure}

In a $n-1$ parameter family of fat graphs, a fat graph $\Gam$
could be a member of two Stasheff associahedra $A_n$. Thus, in
order to have a well-defined framing, we will subdivide the
associahedron making it into a union of a collar neighborhood of
the boundary and a smaller half-sized copy of the associahedron
which we call $A_n^{1/2}$. (See Figure~\ref{MMMW:fig:501}.)

More precisely, $A_n^{1/2}\subseteq A_n=|\A_n|$ is defined to be
the set of all
\[
    (t,\Gam_\ast)\in\Delta^p\times\{\Gam_0\to\cdots\to\Gam_p\}
\]
with the properties:
\begin{enumerate}
  \item[(a)]The last tree $\Gam_p$ contains a vertex of valence $n$ and
  \item[(b)]$t_p\geq\frac12$ ($1$ minus the superscript of
  $A_n$).
\end{enumerate}
Note: An edge of $\Gam_i$ that collapses in $\Gam_{i+1}$ has
length $t_0+\cdots+t_i$ in $\Gam(t)$. Thus $t_p\geq\frac12$ is
equivalent to saying that all edges which collapse in $\Gam_p$
have length $\leq\frac12$ in $\Gam(t)$.

Now we restrict to the case $n=2k+2$. The plan is to take any $2k$
parameter family of fat graphs, subdivide any associahedron
$A_{2k+2}$ which occurs and take the radial framing of the
half-sized subcomplex $A_n^{1/2}$ about the center of mass. It
should be obvious that no fat graph lies in more than one such
subcomplex.

We go through the details. Let $\C_{2k}$ denote the full
subcategory of $\Fat$ consisting of fat graphs of codimension
$\leq2k$ which do not lie in the Witten cycle $W_{k}$, ie, do
not contain vertices of valence $2k+3$. Let $\Gam(t)$,
$t\in|\C_{2k}|$ be the corresponding family of metric fat graphs.

Let $P_{2k+2}^{1/2}$ denote the subset of $|\C_{2k}|$
corresponding to $A_{2k+2}^{1/2}$, ie, the set of all
\[
    (t,\Gam_\ast)\in\Delta^p\times\{\Gam_0\to\cdots\to\Gam_p\}
    \subseteq|\C_{2k}|
\]
with the properties:
\begin{enumerate}
  \item[(a)]$\Gam_p$ contains a vertex of valence $2k+2$ and
  \item[(b)]$t_p\geq\frac12$.
\end{enumerate}
For each $t\in P_{2k+2}^{1/2}$ let $A(t)$ be the tree in $\Gam(t)$
which collapses to the $2k+2$--valent vertex of $\Gam_p$. Then
every edge in $A(t)$ has length $\leq\frac12$. If $\Gam(t)$ has a
massive point then it must be the center of mass of $A(t)$ so
$A(t)\subseteq A_j(t)$ for any $j$ so that $t\in K_j$ (from the
previous subsection). Consequently, Theorem~\ref{WWWM:thm:extended
framed graph theorem} allows us to impose the radial framing on
$A(t)$ centered at the center of mass $c(t)$.

We note that, by definition, the center of mass is either a vertex
or the midpoint of the edge on which it lies. Consequently, when
this edge collapses, both endpoints collapse to $c(t)$
simultaneously. Thus the worse case (\ref{MMMW:eq:worst case no
2}) does not occur since $|C(v)|=2$ implies $v=c(\Gam_0)$ which
implies that $|C(f_1(v))|\geq4$. Therefore, the only contribution
to the $0$--cell cocycle $z^k$ comes from the associahedron
$A_{2k+3}$.

This proves:

\begin{thm}\label{MMMW:thm:Witten conjecture without coefficient}
The $0$--cell cocycle $z^k$ is non-zero only inside the Stasheff
associahedron $A_{2k+3}$ with a standard framing on its boundary.
Consequently, the adjusted Miller--Morita--Mumford class $\wt{\k}_k$
is dual to the Witten cycle $[W_{k}]$ and thus proportional to
$[c_\Fat^k]$ in $H^{2k}(\Fat;\QQ)$.
\end{thm}

%
%

\section{The calculation}

The adjusted Miller--Morita--Mumford class $\wt{\k}_k$ is given by
the $0$--cell cocycle $z^k$ on framed fat graphs which is zero
outside of the Stasheff associahedron $A_{2k+3}=|\A_{2k+3}|$. In
this section we construct an explicit framing for this
associahedron and use it to calculate the proportionality constant
between $\wt{\k}_k$, $[c_\Fat^k]$ and $[W_{k}]^\ast$. Namely:

\begin{thm}\label{Calc:thm:proportionality constant between MMM
and b}
\[
    \wt{\k}_k=-\frac12[c_\Fat^k]
    =(-1)^{k+1}\frac{(k+1)!}{(2k+2)!}[W_{k}]^\ast.
\]
\end{thm}

\subsection{Framing the associahedron $A_{2k+3}$}

First we subdivide $A_{2k+3}$ so that we have a collar
neighborhood $\d A_{2k+3}\times I$ of its boundary $\d A_{2k+3}$
and a half-sized copy $\halfbig$ of the standard associahedron
inside. On this half-sized copy we take the radial framing about
the center of mass. Since $2k+3$ is odd this center of mass will
always be at a vertex. On the outside boundary we take the framing
given in the last section.

\begin{figure}[ht!]
\cl{\includegraphics[width=2.2in]{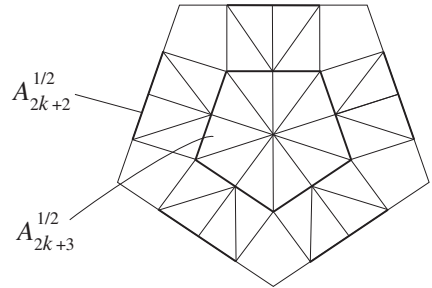}}
\caption{Subdivision of $A_{2k+3}$}\label{Calc:fig:502}
\end{figure}

The outside boundary $\d A_{2k+3}$ contains $2k+3$ copies of the
associahedron $A_{2k+2}$. Following the procedure outlined above
we are required to subdivide each of these associahedra to form a
half-sized copy $\halfsmall$. Then we take the radial framing
centered at the center of mass of the tree which collapses to the
$2k+2$ valent vertex in the center of $\halfsmall$.
Figure~\ref{Calc:fig:502} gives an accurate picture when $k=1$.

\subsection{Value of $z^k$ on $\halfbig$}

We first compute the value of the $0$--cell cocycle on the
half-sized associahedron $\halfbig$. Since we are taking the
radial framing only the center of mass is an unpaired $0$--cell.
Consequently, the only simplices in $\halfbig$ on which the
$0$--cell cocycle $z^k$ has a chance of being non-zero are the ones
in which at each step, $\Gam_i\to\Gam_{i+1}$, one more edge-vertex
pair collapses into the center of mass.

Each time another vertex collapses to the center of mass another
region is added to the middle. (Recall that a \emph{region} is a
gap between consecutive leaves.) There are $2k+3$ regions. We
start with three in the middle $a_1,a_2,a_3$ and we add regions
$b_1,b_2,\cdots,b_{2k}$ one at a time. The positions of these
regions completely determine the shape of the tree in $\Gam_0$
which collapses to the $2k+3$ valent vertex of $\Gam_{2k}$ and the
order of collapse determines (the collapsing part of) every graph
$\Gam_i$.

We will first count the number of $2k$--simplices in $\halfbig$
obtained in this way. Then we will multiply by the average value
of the $0$--cell cocycle on each term.

There are $2k+1$ choices for $a_1$. Fixing $a_1=0$ we have:
\begin{equation}\label{Calc:a1 a2 a3}
  1\leq a_2\leq k+1<a_3\leq a_2+k+1
\end{equation}
since $a_1,a_2,a_3$ are in cyclic order and there are no gaps
greater than $k$ between two consecutive $a$'s. (Otherwise the
center of mass would be in that gap.) Consequently, there are
$k+1$ choices for $a_2$ and the number of choices for $a_3$ is
$a_2$ which has an average value of $\frac{k+2}{2}$. So there are
\[
    \frac13(2k+3)(k+1)\left(\frac{k+2}{2}\right)
\]
choices for $(a_1,a_2,a_3)$. The factor $\frac13$ comes from the
fact that the numbering of the $a$'s is only well-defined up to
cyclic order.

After the $a$'s are chosen, the $b$'s can be chosen arbitrarily.
Thus there are $(2k)!$ choices for $b_1,\cdots,b_{2k}$ making a
total of
\[
    \frac16(2k+3)(k+1)({k+2})(2k)!
\]
terms. The orientation of the $2k$--simplex is given by
(Definition~\ref{FG:def:simplicial orientation of A2k+3})
\[
    \sgn(a_1,a_2,a_3,b_1,\cdots,b_{2k}).
\]
So the value of the $0$--cell cocycle on $\halfbig$ is
\begin{equation}\label{Calc:eq:0-cell cocycle on half big}
    \frac16(2k+3)(k+1)({k+2})(2k)!
    \frac{(-1)^kk!E(Z)}
    {(2k)!(2k+3)!/2}=
    \frac{(-1)^k(k+2)!E(Z)}
    {3(2k+2)!}
\end{equation}
where $E(Z)$ is the expected value of
\[
    Z=\sum_{i=1}^3\sgn(a_i,b_1,\cdots,b_{2k})
    \sgn(a_1,a_2,a_3,b_1,\cdots,b_{2k}).
\]
As in section 3 it is easier to compute
\[
    Z+1=(1+(-1)^P)(1+(-1)^Q)
\]
where $P=a_2-a_1-1$ ($=a_2-1$ when $a_1=0$) and $Q=a_3-a_2-1$. The
random variable $P$ takes values $0\leq p\leq k$ with probability
\[
    P(P=p)=\frac{p+1}
    {(k+1)(k+2)/2}.
\]
The random variable $Q=a_3-a_2-1$ takes values
$
    k-p\leq q\leq k
$
with conditional probability
\[
    P(Q=q|P=p)=\frac1{p+1}.
\]
Thus
\[
    P(Q\text{ is even}|P=p)=
  \begin{cases}
    \frac12 & \text{for $p$ odd}, \\
    \frac12+\frac{(-1)^k}{2(p+1)} & \text{for $p$ even}.
  \end{cases}
\]
So
\[
    E(Z+1|P=p)=
  \begin{cases}
    0 & \text{for $p$ odd}, \\
    2+\frac{(-1)^k2}{p+1} & \text{for $p$ even}
  \end{cases}
\]
and
\begin{align*}
      E(Z+1)&=\sum E(Z+1|P=2j)P(P=2j)\\
      &=\sum_{j=0}^{\left[\frac{k}{2}\right]}
      \left(\frac{2(2j+1+(-1)^k)}{2j+1}\right)
      \left(\frac{2j+1}{(k+1)(k+2)/2}\right)\\
      &=\sum_{j=0}^{\left[\frac{k}{2}\right]}
      \frac{4(2j+1+(-1)^k)}{(k+1)(k+2)}.
\end{align*}

{\bf Case 1}\qua  If $k$ is odd then
\[
    E(Z)=\frac4{(k+1)(k+2)}
    \left(\frac{k-1}{2}\right)
    \left(\frac{k+1}{2}\right)-1=\frac{-3}{k+2}
\]
so
\begin{equation}\label{Calc:eq:value of zk halfbig for k odd}
  z^k(\halfbig)=\frac{(-1)^k(k+2)!}{3(2k+2)!}\cdot\frac{-3}{k+2}
  =(-1)^{k+1}\frac{(k+1)!}{(2k+2)!}.
\end{equation}

{\bf Case 2}\qua If $k$ is even then
\[
    E(Z)=\frac4{(k+1)(k+2)}
    \left(\frac{k+2}{2}\right)
    \left(\frac{k+4}{2}\right)-1=\frac{3}{k+1}
\]
so
\begin{equation}\label{Calc:eq:value of zk halfbig for k even}
  z^k(\halfbig)=\frac{(-1)^k(k+2)!}{3(2k+2)!}\cdot\frac{3}{k+1}
  =(-1)^k\frac{k!(k+2)}{(2k+2)!}.
\end{equation}

\subsection{Value of $z^k$ on the collar}

The only $2k$--simplices in the collar which have a chance of
giving a non-zero value under $z^k$ are the ones in which $\Gam_0$
contains a $0$--cell in the center of an edge $e_0$ which is the
center of mass of an $A_{2k+2}$ tree. The next graph $\Gam_1$ must
be in the boundary $\d\halfbig$ where the center of mass shifts
over to one of the endpoints $v_1$ of $\Gam_1$ which must be equal
to $\Gam_0$.
\begin{figure}[ht!]
\cl{\includegraphics[width=1.4in]{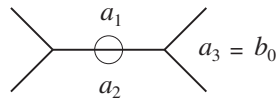}}
\caption{The center of mass (circle) shifts to
$a_3=b_0$.}\label{Calc:fig:503}
\end{figure}

\begin{figure}[ht!]
 \cl{\includegraphics[width=2in]{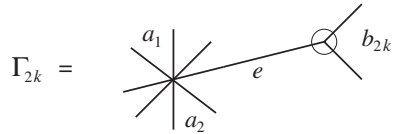}}
\caption{$b_{2k}$ is left out.}\label{Calc:fig:504}
\end{figure}

Using consistent notation with our previous discussion we label
the first three angles $a_1,a_2,a_3$ as shown in
Figure~\ref{Calc:fig:503} where the circle represents the center
of mass for $\halfsmall$ and the right hand vertex is the center
of mass for $\halfbig$. The next regions to converge to the center
of mass are $b_1,\cdots,b_{2k-1}$ with the region $b_{2k}$ being
left out. (See Figure~\ref{Calc:fig:504}.) We note that the
circled vertex should form a collapsing pair with the long edge
$e$.

In order for the center of mass of $\halfsmall$ to be on the edge
separating $a_1,a_2$ we must have $k$ regions on one side and
$k+1$ regions on the other and $b_{2k}$ must be in the part with
$k+1$ regions. The region $a_2=b_0$ must also be on the same side
as $b_{2k}$ and our convention is that $a_1,a_2,a_3$ are in cyclic
order. The other regions $b_1,\cdots,b_{2k-1}$ can be places
arbitrarily. The number of possible configurations is:
\begin{equation}\label{Calc:eq:number of terms in collar}
    (2k+3)(k+1)k(2k-1)!
\end{equation}
($a_1$ can be chosen arbitrarily, there is no choice about $a_2$,
there are $k+1$ choices for $b_0=a_3$ and only $k$ choices for
$b_{2k}$.)

\begin{figure}[ht!]
\cl{\includegraphics[width=3in]{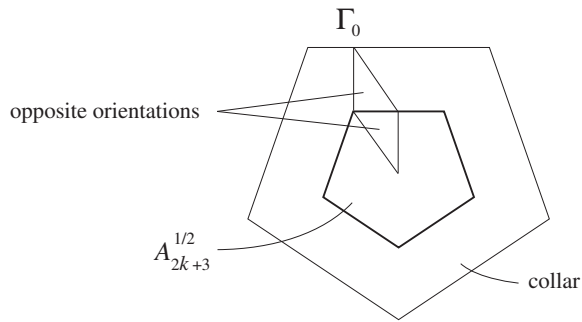}}
\caption{The orientation is reversed in the
collar.}\label{Calc:fig:505}
\end{figure}

The orientation of this $2k$--simplex is given by
\[
    -\sgn(a_1,a_2,b_0,\cdots,b_{2k})
\]
since $\Gam_0$ is on the wrong side of the odd dimensional face
$\Gam_1\to\cdots\to\Gam_{2k-1}$. (See Figure~\ref{Calc:fig:505}.)

The average value of the $0$--cell cocycle $z^k$ on the simplex is
\begin{equation}\label{Calc:eq:zk on simplex in collar}
  \frac{(-1)^kk!}{(2k)!}\frac{E(W)}{(2k+2)!}
\end{equation}
where $E(W)$ is the expected value of
\[
    W=-\sum_{i=1}^2\sgn(a_i,b_0,\cdots,b_{2k-1})\sgn(a_1,a_2,b_0,\cdots,b_{2k}).
\]
If $k$ is odd then $W=0$ since the summands
$\sgn(a_i,b_0,\cdots,b_{2k-1})$ have opposite sign.

If $k$ is even then
\begin{equation}\label{Calc:eq:value of E(W)}
  E(W)=\frac{-2}{k+1}
\end{equation}
since the number of $b$'s between $b_{2k}$ and $a_1$ or $a_2$ is a
random variable $N$ taking the values $0,\cdots,k$ with equal
probability and $
    W=-2(-1)^N.
$

Consequently, when $k$ is even, the total value of $z^k$ on the
collar of $A_{2k+3}$ is given by multiplying (\ref{Calc:eq:number
of terms in collar}) and (\ref{Calc:eq:zk on simplex in collar})
using (\ref{Calc:eq:value of E(W)}):
\[
    -(2k+3)(k+1)k(2k-1)!\frac{(-1)^kk!}{(2k)!}\frac{-2}{k+1}
    =(-1)^{k+1}\frac{k!(2k+3)}{(2k+2)!}.
\]
Adding this to (\ref{Calc:eq:value of zk halfbig for k even})
gives:
\[
    \frac{(-1)^kk!}{(2k+2)!}[k+2-(2k+3)]=(-1)^{k+1}\frac{(k+1)!}{(2k+2)!}
\]
when $k$ is even.

When $k$ is odd we get the same answer (\ref{Calc:eq:value of zk
halfbig for k odd}) since $E(W)=0$. This completes the proof of
Theorem~\ref{Calc:thm:proportionality constant between MMM and b}.

%
%
%

\Addresses\recd
\end{document}